\documentclass[11pt,a4paper,reqno]{amsart}
\usepackage{graphicx} %
\usepackage[table]{xcolor}
\usepackage[normalem]{ulem}

\usepackage{amsmath,amssymb}
\usepackage{xcolor}

\usepackage{mathtools}
\usepackage{pgfplots}
\usepackage{caption}
\usepackage{subcaption}
\pgfplotsset{compat=1.17}
\usepackage{pgfplotstable}
\usepackage{tikz}
\usepackage{multicol}
\usepackage{amsfonts}
\usepackage{multirow}

\usepackage{nicefrac}
\usepackage{etoolbox} 

\usepackage[foot]{amsaddr}

\AtBeginEnvironment{tabular}{%
    \rowcolors{2}{black!10}{}
    }
\usepackage{enumitem}
\allowdisplaybreaks
\pgfplotscreateplotcyclelist{mycyclelist}{
    {blue, solid, mark=o, line width=1pt},
    {red, solid, mark=square, line width=1pt},
    {brown, dashdotted, mark=*, line width=1pt,mark options={solid}},
    {orange, dashdotted, mark=diamond*, line width=1pt,mark options={solid}},
    {green!60!black, dashed, mark=triangle*, line width=1pt,mark options={solid},mark size=3pt},
    {purple, dotted, mark=o, line width=1pt,mark options={solid}},
    {cyan!80!black, dotted, mark=+, line width=1pt,mark options={solid}}
}

\pgfplotscreateplotcyclelist{mycyclelist2}{
    {brown, dashdotted, mark=*, line width=1pt,mark options={solid}},
    {orange, dashdotted, mark=diamond*, line width=1pt,mark options={solid}},
    {green!60!black, dashed, mark=triangle*, line width=1pt,mark options={solid},mark size=3pt},
    {purple, dotted, mark=o, line width=1pt,mark options={solid}},
    {cyan!80!black, dotted, mark=+, line width=1pt,mark options={solid}}
}

\usepackage{url}
\usepackage{algorithm}
 \usepackage{algorithmic}

\newcommand{\R}{\mathbb{R}}
\newcommand{\ddiv}{\operatorname{div}}
\newcommand\norm[1]{\|#1\|}
\newcommand{\calH}{\mathcal{H}}
\newcommand{\calS}{\mathcal{S}}
\newcommand{\calA}{\mathcal{A}}
\newcommand{\calB}{\mathcal{B}}
\newcommand{\calC}{\mathcal{C}}
\newcommand{\calW}{\mathcal{W}}

\newcommand{\calM}{\mathcal{M}}
\newcommand{\calP}{\mathcal{P}}
\newcommand{\dom}{\operatorname{dom}}

\newtheorem{theorem}{Theorem}
\newtheorem{remark}[theorem]{Remark}
\newtheorem{proposition}[theorem]{Proposition}

\newtheorem{lemma}[theorem]{Lemma}

\usepackage{amsmath,amssymb,amsfonts}
\usepackage{latexsym}
\usepackage{mathrsfs}
\usepackage{bbm}
\usepackage{xcolor}
\usepackage{tikz}
\usepackage{url}
\usepackage{stmaryrd}
\usepackage{dsfont}
\usepackage{placeins}
\usepackage{lipsum}
\usepackage{amsfonts}
\usepackage{graphicx}
\usepackage{epstopdf}
\usepackage{amsopn}
\ifpdf
  \DeclareGraphicsExtensions{.eps,.pdf,.png,.jpg}
\else
  \DeclareGraphicsExtensions{.eps}
\fi

\allowdisplaybreaks

\title[Iterative solvers for partial differential equations with dissipative structure]{Iterative solvers for partial differential equations with dissipative structure: Operator preconditioning and Optimal Control}
\author{Volker Mehrmann$^{1}$}
\address{$^{1}$Institute of Mathematics, Technische Universität Berlin, Germany}
\author{Manuel Schaller$^{2}$}
\author{Martin Stoll$^{2}$}
\address{$^{2}$Faculty of Mathematics, Chemnitz University of Technology, Germany}
\thanks{Manuel Schaller is supported by the German Research Foundation (DFG) under project number 519323897.}

\begin{document}

\begin{abstract}
  This work considers the iterative solution of large-scale problems subject to non-symmetric matrices or operators arising in discretizations of (port-)Hamiltonian partial differential equations. We consider problems governed by an operator $\calA=\calH+\calS$ with symmetric part $\calH$ that is positive (semi-)definite and skew-symmetric part $\calS$.  Prior work has shown that the structure and sparsity of the associated linear system enables Krylov subspace solvers such as the generalized minimal residual method (GMRES) or short recurrence variants such as Widlund's or Rapoport's method using the symmetric part $\calH$, or an approximation of it, as preconditioner. In this work, we analyze the resulting condition numbers, which are crucial for fast convergence of these methods, for various partial differential equations (PDEs)  arising in diffusion phenomena, fluid dynamics, and elasticity. We show that preconditioning with the symmetric part leads to a condition number uniform in the mesh size in case of elliptic and parabolic PDEs where $\calH^{-1}\calS$ is a bounded operator. Further, we employ the tailored Krylov subspace methods in optimal control by means of a condensing approach and a constraint preconditioner for the optimality system. We illustrate the results by various large-scale numerical examples and discuss efficient evaluations of the preconditioner, such as incomplete Cholesky factorization or the algebraic multigrid method.\smallskip
    
\end{abstract}
\maketitle

\section{Introduction}\label{sec:intro}
In this work, we investigate iterative solvers for problems governed by linear operators $\calA : \dom(\calA)\supset X\to X$ on a Hilbert space $X$ subject to the {(formal) splitting
\begin{equation}\label{eq:splitting}
    \calA = \calH + \calS %
\end{equation}
where $\calS:X \supset \dom(\calS)\to X$ is skew-symmetric (or skew-adjoint) and $\calH : X\supset \dom(\calH) \to X$ is symmetric (or self-adjoint) and nonnegative.} %
This structure is ubiquitous in stationary or time-dependent partial differential equations (PDEs) governed by dissipative operators such as port-Hamiltonian systems. Therein, the positive definiteness of the symmetric part $\mathcal{H}$ corresponds to coercivity of a governing even-order differential operator such as in diffusion problems, fluid dynamics or elasticity with strong damping.

Previous work~\cite{diab2023flexible,GuduLies22,ManM19,ManM21,MangMehr25}, has suggested to leverage a skew-sym\-met\-ric/sym\-met\-ric splitting~\eqref{eq:splitting} for efficient solutions of matrix problems, i.e., via left-preconditioning in the generalized minimal residual method (GMRES), see e.g. \cite{saad1986gmres}, using the symmetric part, or via Widlund's or Rapoport's method, \cite{rapoport1978nonlinear,widlund1978lanczos}, which constitute short-term recurrence Krylov subspace methods. As usual for iterative methods, the convergence rate strongly depends on the condition number, respectively the spectral width, of the preconditioned system, see e.g., \cite{GuduLies22,szyld1993variational}. In the case of an underlying infinite-dimensional problem (e.g., a partial differential equation) governed by an operator \eqref{eq:splitting} with invertible symmetric part $\calH$, this condition number is determined by Galerkin projections of 
\begin{equation*}
    \mathcal{H}^{-1}\mathcal{A} = I + \mathcal{H}^{-1}\mathcal{S}
\end{equation*}
onto a suitable finite element space. In this work, we will analyze this approach through the lens of operator preconditioning and provide applications in large-scale optimal control of PDEs.

Operator preconditioning is a powerful framework to assess the conditioning of Galerkin projections by means of the analysis of the underlying partial differential equation in function space. There is a rich literature for this topic and we refer to the overview articles~\cite{hiptmair2006operator,axelsson2009equivalent} focusing on elliptic problems, the work \cite{mardal2011preconditioning} for an abstract and operator-theoretic approach to problems in Hilbert spaces. In \cite{gunnel2014note}, the equivalence of choosing a preconditioner and the choice of a suitable inner product for conjugate gradient and minimal residual methods is presented. For an operator preconditioning approach to equality-constrained optimization, we refer to \cite{schiela2014operator}. Preconditioning with the symmetric part for Navier-Stokes equations was analyzed in \cite{benzi2011relaxed}. %
For a general overview of preconditioning, in particular in the context of partial differential equations, we refer to the overview article~\cite{wathen2015preconditioning} (see in particular Section 6.4 for preconditioning with the symmetric part), the book~\cite{elman2014finite}, the overview article~\cite{malek2014preconditioning} focusing particularly on the conjugate gradient method or the book~\cite{toselli2004domain} for domain decomposition-based preconditioners. 

In this work, we leverage the structure of the underlying operator equation and its discretized version in particular for preconditioning in simulation and optimal control of partial differential equations. In this context, the accretivity of the underlying equation (or equivalently, the dissipativity) plays a crucial role and is, in a wide range of applications induced by a dissipative Hamiltonian or port-Hamiltonian structure, see \cite{AchAM23a,AchAMN24,mehl2022matrix,MehMW25,ZwaM24} for an
analysis of the spectral properties of matrix and operator pencils with dissipative Hamiltonian structure. We also briefly mention other works, where a dissipative Hamiltonian structure is used in (numerical) linear algebra. In \cite{diab2023flexible,GuduLies22,ManM19,ManM21} short-term recurrence Krylov subspace methods for matrices occurring in dissipative Hamiltonian problems are discussed, see also the recent work \cite{MangMehr25} for a treatment of systems with saddle point structure. In~\cite{bartel2025splitting,bartel2023operator}, splitting and dynamic iteration methods leveraging the modular port-Hamiltonian structure were suggested, where the dissipativity particularly provides a monotonous bound on the convergence. Peaceman-Rachford-type splitting methods in function space were suggested in \cite{farkas2023operator} for simulation and in~\cite{farkas2025dissipativity} for optimal control.

\textbf{Contribution and outline.} This paper features two main contributions, the first being an analysis of particular PDEs through the lens of operator preconditioning using the symmetric part, the second focusing on structure-exploiting numerical methods for the solution of the linear systems arising in optimal control. In the first part in Section~\ref{sec:conditioning}, we investigate the conditioning for a wide range of example problems such as advection-diffusion-type problems, fluid dynamics problems such as Stokes and Oseen equations or problems from mechanics such as the wave equation and a beam equation. We will observe that depending on the functional analytic structure of the equation, in particular if $\mathcal{H}^{-1}\mathcal{S}$ is a bounded operator, the preconditioning leads to a condition number that is uniform in the mesh size. As a second contribution in Section~\ref{sec:optcont}, we illustrate how iterative solvers (such as GMRES preconditioned with the symmetric part, Widlund's or Rapoport's method) may be implemented in numerical optimal control of PDEs. Particularly, as the adjoint equation inherits the structure of the state equation (e.g, accretive, dissipative or port-Hamiltonian), we leverage the structure-exploiting methods for solving state and adjoint equation. Using this, we first suggest an elimination to the control, in which the corresponding optimality system is solved with a conjugate gradient method, and in which the inner evaluation of the optimality system is performed with state and adjoint equation solves including the above mentioned methods. Second, we will use the methods in a constraint preconditioner in a conjugate gradient method applied to the full optimality system. By proposing efficient approximations to $\mathcal{H}^{-1},$ such as incomplete Cholesky factorizations or multigrid methods, we show that the proposed preconditioning technique may be efficiently implemented within the optimal control problem. %

\textbf{Notation.} We adopt standard notations from operator theory (see~\cite{EngelNagel2000}), functional analysis, and partial differential equations (see~\cite{adams2003sobolev}). Let $(X,\langle\cdot,\cdot \rangle)$ be a Hilbert space. We say that an operator $\calA : \dom(A) \subset X\to X$ is \emph{accretive} if $\langle \calA x,x\rangle \geq 0$ for all $x\in \dom(\calA)$, and \emph{strictly accretive} if $\langle \calA x,x\rangle \geq \mu \|x\|^2$ for some $\mu > 0$ and all $x\in \dom(\calA)$. We call $\calA$ \emph{(strictly) dissipative}, if $-\calA$ is (strictly) accretive. If $\dom(\calA)$ is dense in $X$, the \emph{adjoint of $\calA$} is defined by $\calA^*\colon \dom(A^*)\subset X \to X$ with
	\begin{equation}\label{eq:domAstern}
	     \dom(\calA^*) \coloneqq \left\{y \in X \mid \mbox{\rm there is } z\in X: \langle y,\calA x \rangle_{X} = \langle z,x\rangle_{X}\  \mbox{\rm for all } x \in \dom(A)\right\}.
	\end{equation}
Due to the density of $\dom(\calA)$ in $X$, the element $z$ in this set is uniquely determined and we set $\calA^* y \coloneqq z$. We say that $\calA$ is \emph{symmetric (resp.\ skew-symmetric)}, denoted by $\calA\subset \calA^*$ (resp.\ $\calA \subset-\calA^*)$ if or all $x\in \dom(\calA)$, $\calA x=\calA ^*x$ (resp.\ $\calA x=-\calA^*x$). We say that $\calA$ is \emph{self-adjoint (resp.\ skew-adjoint)}, denoted by $\calA=\calA^*$ (resp.\ $\calA = - \calA^*$) if $\calA$ is symmetric (resp.\ skew-symmetric) and $\dom(\calA) = \dom(\calA^*)$.

\medskip

\textbf{Code availability and reproducibility.} The code for all algorithms and numerical experiments conducted in this work, as well as the involved matrices and the files to generate the discretizations via the finite element library FEniCS~\cite{alnaes2015fenics} is provided in the repository 
\begin{center}
    \url{https://github.com/maschaller/indefinite_solvers}.
\end{center}

\medskip

Let $\Omega \subset \R^d$. For a function \( f : \Omega \to \mathbb{R} \), we denote its gradient by \( \nabla f \). For a vector field \( g : \Omega \subset \mathbb{R}^d \to \mathbb{R}^d \), we denote its divergence by $\operatorname{div} g = \sum_{i=1}^d \frac{\partial g_i}{\partial \omega_i}.$
We denote by $L^2(\Omega;\R^d)$ resp.\ $L^\infty(\Omega;\R^d)$) the space of square-integrable resp.\ measurable and essentially bounded (equivalence classes of) functions $f:\Omega\to \R^d$. By \( H^1(\Omega;\R^d) \) we denote the standard Sobolev space of (equivalence classes of) square integrable functions $f:\Omega\to \R^d$ with square-integrable first derivatives and the subspace \( H^1_0(\Omega;\R^d) \) consists of functions in \( H^1(\Omega;\R^d) \) with vanishing trace on the boundary.  In case $d=1$, we abbreviate by omitting the second argument, e.g., $L^2(\Omega)= L^2(\Omega;\R)$.
Moreover, \( H(\operatorname{div}, \Omega) \) denotes the space of vector fields 
\( f : \Omega \subset \mathbb{R}^d \to \mathbb{R}^d \) with components in \( L^2(\Omega) \) 
whose divergence (in the weak sense) also belongs to \( L^2(\Omega) \).

\section{Specialized solvers leveraging symmetric/skew-symmetric splittings}\label{sec:solvers}
In the following, we will briefly recall the solvers discussed in \cite{diab2023flexible,GuduLies22} for a linear system
\begin{equation}\label{eq:discsys}
    (H+S)x = b,
\end{equation}
in which $H\in \R^{n\times n}, S\in \R^{n\times n}$ typically result as Galerkin projections of operators $\mathcal{H}$ and $\mathcal{S}$ from \eqref{eq:splitting}, respectively. In particular, we assume in the following that $H=H^\top >0$ and $S=-S^\top$. The case of possibly non-invertible $H=H^\top \geq 0$ will be discussed later by a Schur complement approach.

To design iterative methods for the finite dimensional non-symmetric problem~\eqref{eq:discsys}, one could choose a very general approach such as \textsc{GMRES} \cite{saad1986gmres} endowed with a suitable preconditioner. The downside of \textsc{GMRES} and the motivation for its variants is that new vectors need to be orthogonalized against all basis vectors within the Krylov subspace. For purely symmetric or skew-symmetric system matrices one can rely on short-term recurrence methods such as \textsc{MINRES} or \textsc{CG} that only require orthogonalization against a small number of basis vectors. We now briefly introduce methods that share this favorable property and rely on the the fact that the system matrix is symmetric or skew symmetric  in a non-standard inner product \cite{faber2008faber,stoll2008combination}. In our case, the matrix $H^{-1}S$ is skew-symmetric in the inner product defined by the symmetric part $H$, that is,
\begin{align}\label{eq:symmetric}
\begin{split}
\langle H^{-1}Sx,x\rangle_H&=x^\top(H^{-1}S)^\top Hx=x^\top S^\top H^{-1}Hx\\ %
&=x^\top HH^{-1}S^\top x=-x^\top HH^{-1}Sx=-\langle x,H^{-1}Sx\rangle_H.
\end{split}
\end{align}
Methods that are tailored to exploit this skew-symmetry for a short-term recurrence were introduced by Widlund \cite{widlund1978lanczos} and Rapoport \cite{rapoport1978nonlinear}. The main idea is to built up a Krylov subspace via a Lanczos relation
\[
H^{-1}SV_{k}=V_{k+1}T_{k+1,k},
\]
where $V_k$ and $V_{k+1}$ store the basis vectors for the Krylov subspace and $T_{k+1,k}\in\R^{k+1,k}$ is a tridiagonal coefficient matrix. Here, in view of \eqref{eq:symmetric}, the Lanczos method is based on orthogonality in the inner product defined by $H$. The methods by Widlund and Rapoport are then obtained by minimizing a desired quantity such as the residual $r_k$ or the error $x-x_k$ in a norm induced by $H$.  In more detail, setting $\hat b = H^{-1} b$, Widlund's method corresponds to solving
\[
(I_k+T_{k,k})y_k=\norm{\hat{b}}_He_1. 
\]
which yields the iterate $x_k=V_ky_k.$ For Rapoport's method a residual minimization in the $H^{-1}$ norm results in the solution of the $k\times k$ system
\[
T_{k+1,k}^{\top}T_{k+1,k}y_k=\norm{\hat{b}}_HT_{k+1,k}^\top e_1,
\]
where again $x_k=V_ky_k$. We note that similar search spaces are used when applying \textsc{GMRES} to the left-preconditioned system $(I+H^{-1}S)x=H^{-1}b$. GMRES, however, requires an orthogonalization against \emph{all} previous basis vectors in the Krylov space and thus can become quite slow if the preconditioner is not very effective. To this end, restarted variants of \textsc{GMRES} or methods with smaller memory such as Loose \textsc{GMRES} (LGMRES) have been proposed~\cite{baker2005technique}. %
If the preconditioner is evaluated inexactly, e.g., when using only a few steps of a multigrid method, flexible solvers provide a framework to ensure convergence and low-recurrence in Lanczos methods despite inexactness. We refer the reader to the work~\cite{diab2023flexible} which suggests various flexible variants for preconditioning the system \eqref{eq:discsys} with the symmetric part. 

All the discussed methods, i.e., Widlund's and Rapoport's method, as well as GMRES, are iterative Krylov subspace methods. As such, their convergence strongly depends on the condition number, respectively the spectral width, of the governing matrix $I+H^{-1}S$ as illustrated in \cite[Section 5]{GuduLies22}. To this end, let $\lambda \in \R$ be such that we have the spectral inclusion $\sigma(H^{-1}S)\subset i[-\lambda,\lambda]$ (note that $H^{-1}S$ is $H$-normal due to~\eqref{eq:symmetric} and thus $H^{-1}S$ has imaginary eigenvalues). For Widlund's method, the rate is given by
\begin{equation}\label{eq:wid}
    \frac{\|x - x_{2k}\|_H}{\|x\|_H} \leq 2\left(\frac{\sqrt{1+\lambda^2} - 1}{\sqrt{1+\lambda^2} +1} \right)^k.
\end{equation}
For Rapoport's method, we have convergence of the residual of the form
\begin{equation}\label{eq:rap}
    \frac{\|b-(H+S)x_k\|_{H^{-1}}}{\|b\|_{H^{-1}}} \leq 2\left(\frac{\lambda}{\sqrt{1+\lambda^2}+1}\right)^k.
\end{equation}
In both estimates, the convergence speed strongly depends on the spectral width of $H^{-1}S$, that is, if $\lambda\to \infty$, the convergence deteriorates. For discretizations of partial differential equations, it is thus desirable that this spectral width is bounded uniformly in the mesh, which is closely associated to boundedness of the underlying (preconditioned) operator $\calH^{-1}\calS$. In the following section, we will inspect this boundedness property and associated conditioning for various problem classes.

\section{Spectral properties in Hermitian preconditioning of partial differential equations}\label{sec:conditioning}
In this section, we briefly recall the idea of operator preconditioning in Krylov subspace methods and refer to %
\cite{axelsson2009equivalent,gunnel2014note,hiptmair2006operator,mardal2011preconditioning} for more details. Let $\mathcal{A}:X \supset \dom(\calA) \to X$ be a closed and densely defined but unbounded linear operator. 
Here we will call a densely defined operator $\calA:X\supset \dom(\calA)\to X$ \emph{unbounded} on $(X,\|\cdot\|)$ if there is no $c\geq 0$ such that $\|\calA x\| \leq c\|x\|$ for all $x\in \dom(\calA)$.  

Then, for given $b\in X$, we aim to find $x\in \dom(\calA)$ such that 
\begin{equation}\label{eq:operatorequation}
    \mathcal{A}x = b.
\end{equation}
Since $\calA$ is unbounded on $X$, the Krylov subspace $\mathcal{K}_{m}(\calA,b) := \{b, \calA b, \ldots, \calA^{m-1}b\}$ is not well-defined, as applying powers of $\calA$ is not feasible due to $\calA^2b = \calA \calA b$ but, in general $\calA b \notin \dom(\calA)$.  To alleviate this, in operator (left) preconditioning one introduces a mapping $\calP \in L(X,\dom(\calA))$ and considers the problem
\begin{equation*}
    \calP \calA x = \calP x
\end{equation*}
such that $\calP \calA \in L(\dom(\calA),\dom(\calA))$. We note that $\calA$ is closed, we may render $\dom(\calA)$ a Banach space when endowed with the graph norm $\|x\|_{\dom(\calA)} := \|x\|+\|Ax\|$.
Consequently, the Krylov subspace $\mathcal{K}_{m}(\calP\calA,\calP b) = \{\calP b, \calP \calA \calP b, \ldots, (\calP \calA)^{m-1}\calP b\}\subset \dom(\calA)$ is well-defined.

A simple illustration of operator preconditioning may be obtained by the second-order operator governing the advection-diffusion-reaction equation on $\Omega \subset \R^d$, which we will discuss in detail and in a more general setting in Subsection~\ref{subsec:advec}.  Consider 
the PDE
\begin{equation}\label{eq:adv_diff_stat}
    \mathcal{A}x = -\Delta x +  \mathbf {b}^\top \nabla x + \mathbf {c} x=f
\end{equation}
for $\mathbf {c} >0$ and $\mathbf {b}:\Omega \to \R^d$. Endowing the system with homogeneous Dirichlet boundary conditions, we may define the operator $\calA$ on $X = H^1_0(\Omega)^*$ with $\dom(\calA)=H^1_0(\Omega)$. Here, a suitable preconditioner is $(-\Delta+{\mathbf c} I)^{-1} : X \to \dom(\calA)$ which corresponds to the Riesz isomorphism in $H^1(\Omega)$ \cite{mardal2011preconditioning}. This may be understood as choosing a suitable scalar product for the computation of, e.g., gradients in a (conjugate) gradient method~\cite{gunnel2014note}. A similar argument applies when considering $X=L^2(\Omega)$ with domain $\dom(\calA)=H^1_0(\Omega) \cap H^2(\Omega)$. %

Besides being crucial in well-definedness of the Krylov subspace for the governing infinite-dimensional operator equation \eqref{eq:operatorequation}, preconditioning is a central aspect in ensuring fast convergence of iterative methods. In conjugate gradients for symmetric problems, the convergence is strongly depending on the condition number, %
while in Rapoport's and Widlund's method, we observe a dependence on the spectral width, see \eqref{eq:wid} and \eqref{eq:rap}. If the underlying problem is subject to an unbounded operator, which hence has unbounded spectrum, this implies that Galerkin projections onto suitable spaces, such as finite element spaces $V_{h,k}$ of order $k\in \mathbb{N}$ and mesh size $h>0$, have an inherently bad conditioning in particular when refining the mesh. This may be seen easily with the Cauchy interlacing theorem (or min-max theorem) \cite[Theorem 9.12]{boffi2010finite} in Babu\v{s}ka--Osborn theory, which states that (assuming real eigenvalues)
\begin{equation}\label{eq:interlacing}
    \lambda_i \leq \mu_i \leq \lambda_i + C(i)\mathcal{O}(h^k) \qquad \mbox{\rm for all } 1\leq i \leq N:=\dim V_{h,k},
\end{equation}
where the $\lambda_i$ are the eigenvalues of $\calA$ and the $\mu_i$ are the eigenvalues of the Galerkin projection $P_{V_{h,k}}\calA P_{V_{h,k}}$, both sorted in ascending order, and where $P_{V_{h,k}}$ is the orthogonal projection onto the finite element subspace $V_{h,k}$. If $\calA$ is unbounded, then the spectrum of $\calA$ is, as well, such that $|\lambda_i|\to \infty$ for $i\to \infty$.
Consequently, considering the condition number corresponding to the spectral norm $\kappa_2(P_{V_{h,k}}\calA P_{V_{h,k}}) = \nicefrac{\mu_N}{\mu_1}$ diverges as $|\mu_N| %
\to \infty$ as $h\to 0$ due to \eqref{eq:interlacing} and the unboundedness of $\calA$. %
This divergence may be made precise when considering particular settings. For Galerkin projections of second-order operators such as~\eqref{eq:adv_diff_stat} and piecewise affine linear finite elements, the condition number is of order $h^{-2}$, see \cite[Theorem 1.32]{elman2014finite}. In the numerical examples presented later in this work, we will observe also different orders of the condition number, closely related to the order of the underlying differential operator $\calA$ and its preconditioned variant $\calH^{-1}\calA$.

In the following subsections, we will inspect the suitability of self-adjoint preconditioning for problems governed by operators that admit a self-adjoint/skew-adjoint splitting as in \eqref{eq:splitting}.
In view of the discussion above, to ensure that the self-adjoint part is a suitable preconditioner that in particular renders the condition number, respectively the spectral width, uniformly bounded in terms of the mesh width $h>0$, we have to ensure that
\begin{equation}\label{eq:bdd}
    \calH^{-1} \calA = I + \calH^{-1}\calS \in L(X,X).
\end{equation}
We will observe that this boundedness condition is usually the case for parabolic or elliptic problems as illustrated in Subsection~\ref{subsec:advec} for an advection-diffusion-type problem, in Subsection~\ref{subsec:stokes} for the Stokes equation, in Subsection~\ref{subsec:oseen} for an Oseen equation or in Subsection~\ref{subsec:kv} for a beam equation with structural Kelvin-Voigt damping. 

However, for the wave equation with momentum damping, where the governing derivative operator of highest order is skew-symmetric, we will observe in Subsection~\ref{subsec:wave} that Hermitian preconditioning is not suitable if one aims for a condition number bounded uniformly in the mesh size. 

For the sake of clarity of presentation, we present results for stationary problems. However, we stress that all considerations also carry over to semi-discretizations of dissipative Hamiltonian problems of the form
\begin{equation}\label{eq:timedep}
    \dot  x= -\calM x,
\end{equation}
where $\calM:X\supset\dom(M) \to X$ is an accretive operator, i.e., $\operatorname{Re}\langle x,\calM x\rangle \geq 0$ for all $x\in \dom(\calM)$ which we may formally split via $\calM = \calH_\calM + \calS_\calM$ with symmetric and accretive operator $\calH$ and skew-symmetric operator $\calS$. Note that for general operators, such a splitting can be very intricate, in particular the naive definition $\calH = \tfrac12(\calM + \calM^*)$ may lead to a trivial domain~\cite{arlinskiui2020everything}. However, for particular applications in PDEs, we will make this splitting mathematically precise in the subsequent subsections. Semi-discretization in time, e.g., by the implicit midpoint rule with step size $\delta t>0$, leads to the implicit iteration
\begin{equation}\label{eq:IMP}
    (I+\tfrac{\delta t}{2}\calM)x^+ = (I-\tfrac{\delta t}{2}\calM)x.
\end{equation}
Due to the accretivity of $\calM$, the operator on the left-hand side of~\eqref{eq:IMP} is strictly accretive, i.e., $\langle x,I+\tfrac{\delta t}{2} \calM x\rangle \geq  \|x\|^2$ and (at least formally) may be split into the symmetric and accretive operator $\calH = I + \frac{\delta t}{2}\calH_\calM$ and the skew-symmetric operator $\calS = \frac{\delta t}{2}\calH_\calS$, where $\calH_\calM$ and $\calH_\calS$ are the symmetric resp.\ skew-symmetric parts of $\calM$. Thus, the governing operator is simply a scaled variant of the accretive operator $\calM$ shifted by the identity, which of course does not influence its boundedness and thus neither the conditioning nor the spectral width of Galerkin projections.

We also briefly note that in the remainder of this work, we will sometimes show  the (skew-)symmetry of operators (in contrast to the stronger property of skew-adjointness and self-adjointness). However, when considering Galerkin methods, one usually chooses the subspace $V_h = \operatorname{span}\{\varphi_1,\ldots,\varphi_N\} \subset \dom(\calA)$ and performs a projection of the operator by defining $A\in \R^{N\times N}$ by
\begin{align*}
    A_{ij} := \langle \calA \varphi_i, \varphi_j\rangle_{L^2(\Omega)} \qquad \mbox{\rm for all } i,j=1,\ldots,N
\end{align*}
hence only testing against elements in $\dom(\calA)$. Consequently, skew-symmetry $\calA x = -\calA^* x$ for all $x\in \dom(\calA)$ or symmetry $\calA x = \calA^* x$ for all $x\in \dom(\calA)$ (hence only $\dom(\calA^*) \subset \dom(\calA)$ and not necessarily $\dom(\calA^*) = \dom(\calA)$) is enough to ensure (skew)-symmetry of the resulting matrix. A different domain of the adjoint can be included using Petrov-Galerkin methods, in which the test space differs from the Ansatz space.

\subsection{Stationary advection-diffusion-reaction equation}\label{subsec:advec}
As a first model, we consider a stationary advection-diffusion-reaction equation with homogeneous Dirichlet boundary conditions given by the PDE
\begin{equation}\label{eq:operatoradv}
\begin{split}
\begin{aligned}
  -\ddiv(\nu \nabla x) + \mathbf {b}^{\top}\nabla x + \mathbf{c}x & = f & \quad & \text{on } \Omega, \\
  x            & = 0 & \quad & \text{on } \partial\Omega,
\end{aligned}
\end{split}
\end{equation}
where $x:\Omega \to \R$ models the scalar-valued concentration of a species on a spatial domain $\Omega\subset \R^d$, $d\in \mathbb{N}$ and $f:\Omega \to \R$ is a given source term. The \emph{diffusivity} $\nu\in L^\infty(\Omega,\R^{d\times d})$ is a pointwise symmetric and uniformly positive matrix-valued mapping, i.e., $\nu(\omega)^\top = \nu(\omega)$ for all $\omega\in \Omega$ and there is $\underline \nu > 0$ such that $v^\top \nu(\omega) v \geq \underline \nu \|v\|_{\R^d}^2$ for all $v\in \R^d$ and $\omega \in \Omega$. To ensure accretivity of the governing operator, we assume that $\mathbf {b}\in L^\infty(\Omega;\R^d)$ has a vanishing weak divergence, i.e., $\mathbf{b}\in H(\ddiv,\Omega)$, $\ddiv \mathbf {b} = 0$ a.e.\ on $\Omega$ and that $\mathbf {c}\in L^\infty(\Omega;\R)$ is pointwise nonnegative, i.e., $\mathbf{c}(\omega)\geq 0$ for a.e.\ $\omega \in \Omega$.

In the following Lemma we analyze the structure of the associated operator.
\begin{lemma}\label{lem:struc_advdiff}
    Let $D=\{x\in H^1_0(\Omega)\,|\, \nu \nabla x\in H(\ddiv,\Omega)\}$ and define
    \[
    \begin{aligned}
        \calA &:L^2(\Omega) \supset D \to L^2(\Omega), &\hspace{.3cm}& \calA x = -\ddiv(\nu \nabla x) + \mathbf {b}^\top \nabla x + \mathbf {c}x,\\
        \calH &:L^2(\Omega) \supset D \to L^2(\Omega), &\hspace{.3cm}& \calH x=- \ddiv(\nu \nabla x) + \mathbf {c}x,\\ 
        \mathrm{and}\quad \calS&: L^2(\Omega)\supset D\to L^2(\Omega), &\hspace{.3cm}& \calS x=  \mathbf {b}^\top \nabla x.
    \end{aligned}
    \]
     Then $\calH$ is boundedly invertible and can be split as
    \begin{equation*}
        \calA = \calH + \calS \quad \mathrm{with}\quad \calH=\calH^*, \quad \mathrm{and}\quad \calS \subset -\calS^*,
    \end{equation*}
    where we use the notation as introduced after \eqref{eq:domAstern}.
\end{lemma}

\begin{proof}
We first prove the bounded invertibility of $\mathcal{H}$. First, for all $x\in D$,
\begin{align*}
    \langle \calH x,x\rangle_{L^2(\Omega)} &= -\langle \ddiv(\nu \nabla x),x\rangle_{L^2(\Omega)} + \langle \textbf{c}x,x\rangle_{L^2(\Omega)}\\
    &= \langle \nu \nabla x,\nabla x\rangle_{L^2(\Omega;\R^d)} + \langle \textbf{c}x,x\rangle_{L^2(\Omega)}  %
    \\ 
    &\geq \tfrac12 \min\{\alpha ,1\}\underline{\nu}\|x\|^2_{H^1(\Omega;\R^d)}
\end{align*}
where $\alpha>0$ is the Poincaré constant, i.e., $\|\nabla x\|^2_{L^2(\Omega)} \geq \alpha \|x\|_{L^2(\Omega)}^2$ for all $x\in D\subset H^1_0(\Omega)$. Hence,  by the Lax-Milgram theorem and maximal elliptic regularity~\cite[Section 3]{grisvard2011elliptic}, $\calH$ is boundedly invertible.
We proceed to show self-adjointness of $\calH$. Let $x_1,x_2\in D$. Then, due to the symmetry of $\nu:\R^{d}\to \R^{d\times d}$ and by the previous computations, we have
\begin{align*}
    \langle \calH x_1,x_2\rangle_{L^2(\Omega)}
    &= \langle \nu \nabla x_1, \nabla x_2\rangle_{L^2(\Omega)} + \langle  \textbf{c}x_1,x_2\rangle_{L^2(\Omega)}\\
    &= \langle  \nabla x_1,\nu \nabla x_2\rangle_{L^2(\Omega)} + \langle x_1,  \textbf{c}x_2\rangle_{L^2(\Omega)}\\
    &=-\langle \ddiv(\nu \nabla x_2),x_1\rangle_{L^2(\Omega)} + \langle  \textbf{c}x_2,x_1\rangle_{L^2(\Omega)}= \langle \calH x_2,x_1\rangle_{L^2(\Omega)}
\end{align*}
such that $\calH$ is symmetric. As $\calH$ is boundedly invertible, it is also self-adjoint. This follows from $\calH = (\calH^{-1})^{-1} = (\calH^{-*})^{-1} = \calH^*$, as the bounded inverse of a symmetric operator is always  self-adjoint.

We now compute the adjoint of $\calS$. Let $x_1,x_2\in D$. Then
\begin{align}\label{eq:compS}
    \langle \calS x_1,x_2\rangle_{L^2\Omega)} = \langle \textbf{b}^\top \nabla x_1,x_2\rangle_{L^2(\Omega)} = -\langle x_1,\ddiv(\textbf{b} x_2)\rangle_{L^2(\Omega)} = \langle x_1,\calS^* x_2\rangle_{L^2(\Omega)}
\end{align}
by the definition of the adjoint. This implies that $\dom(\calS^*) = \{x\in L^2(\Omega)\,|\, \textbf{b}x\in H(\ddiv,\Omega)$\}. Hence $\dom(\calS^*)\supsetneq \dom(\calS) = D$. For $x\in \dom(\calS)=D\subset H^1_0(\Omega)$, we may proceed and use the chain rule for the weak divergence and obtain
\begin{align*}
    \calS^*x = -\ddiv(\textbf{b}x) = -(\ddiv \textbf{b}) x  - \textbf{b}^\top \nabla x = -\textbf{b}^\top \nabla x
\end{align*}
as $\textbf{b}$ is divergence-free. Together with \eqref{eq:compS} this implies that $\calS$ is skew-symmetric.
\end{proof}

We note that in the previous result, we could also define $\calS$ on the space $H^1(\Omega)$, which would be the canonical extension as $\calS$ involves only a gradient term.
\begin{proposition}[Preconditioning]\label{prop:adv_diff}
    The preconditioned operator $I + \calH^{-1}\calS: H^1_0(\Omega) \to H^1_0(\Omega)$ is bounded.
\end{proposition}
\begin{proof}
As $\calS : H^1_0(\Omega) \to L^2(\Omega)$ is linear and bounded, and $\calH:D\to L^2(\Omega)$ is boundedly invertible, $\calH^{-1}\calS : H^1_0(\Omega)\to \dom(\calH) =\{x\in H^1_0(\Omega)\,|\, \nu \nabla x\in H(\ddiv,\Omega)\} \hookrightarrow H^1_0(\Omega)$ continuously such that $\calH^{-1}\calS\in L(H^1_0(\Omega),H^1_0(\Omega))$. This yields the assertion.
\end{proof}

In the following remarks we  briefly discuss compactness and weak formulations. %
\begin{remark}\label{rem:weak} 
1)    Note that if the boundary is smooth enough and $\nu \in C^1({\Omega})\cap C(\bar{\Omega})$, then $\dom(\calH) =\{x\in H^1_0(\Omega)\,|\, \nu \nabla x\in H(\ddiv,\Omega)\} = H^2(\Omega)\cap H^1_0(\Omega)$ such that the embedding $\dom(\calH) \hookrightarrow H^1_0(\Omega)$ is compact. Consequently, $\mathcal{H}^{-1}\mathcal{S}$ is also compact. A similar argument may also be employed for less regular boundaries and coefficients, where we have an embedding $\dom(\calH)\hookrightarrow H^{1}_0(\Omega) \cap H^{1+\varepsilon}(\Omega)$ for some $\varepsilon>0$ which is sufficient for compactness. We refer the reader to~\cite[Section 2]{grisvard2011elliptic}. However, we note that this compactness does not carry over to the operator $I+\calH^{-1}\calS$ which is the governing operator of the problem, see \eqref{eq:bdd}.

    2) The operator $\calS$ admits a unique extension $\calS_\mathrm{e}:H^{-1}(\Omega)\supset \dom(\calS_\mathrm{e}) = L^2(\Omega) \to H^{-1}(\Omega)$ defined by
\begin{equation*}
	\langle \calS_\mathrm{e}x,y\rangle_{H^{-1}(\Omega),H^1(\Omega)} = -\langle x,\ddiv (\mathbf{b}x)\rangle_{L^2(\Omega)}
\end{equation*}
for $x\in \dom(\calS_\mathrm{e})$ and $y\in H^1_0(\Omega)$.
The same applies to $\calH$ that admits a unique extension $\calH_\mathrm{e}:H^{-1}(\Omega)\supset \dom(\calH_{\mathrm{e}})= H^1_0(\Omega)$ defined by
\begin{equation*}
\langle \calH x,y\rangle_{H^{-1}(\Omega),H^1(\Omega)} = \langle \nu \nabla x,\nabla y\rangle_{L^2(\Omega)} + \langle \mathbf{c}x,y\rangle_{L^2(\Omega)}
\end{equation*} 
for $x\in \dom(\calH_\mathrm{e})$ and $y\in H^1_0(\Omega)$. These extensions correspond to the \emph{weak form} of the PDE \eqref{eq:operatoradv}. By straightforward adaption of the proof of Proposition~\ref{prop:adv_diff}, we may analogously obtain boundedness of $\calH_\mathrm{e}^{-1}\calS_\mathrm{e}: L^2(\Omega)\to L^2(\Omega)$.
\end{remark}

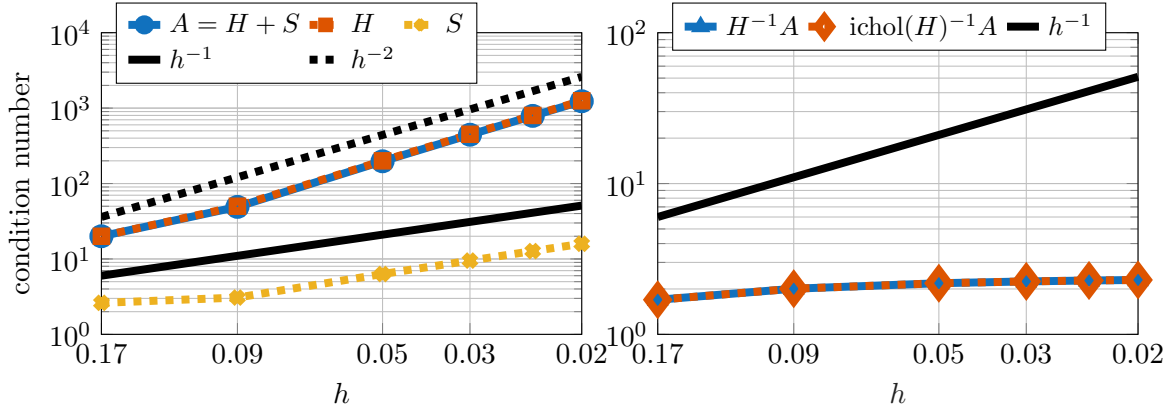
\begin{figure}[htb]
    \centering
\definecolor{mycolor1}{rgb}{0.06600,0.44300,0.74500}%
\definecolor{mycolor2}{rgb}{0.86600,0.32900,0.00000}%
\definecolor{mycolor3}{rgb}{0.92900,0.69400,0.12500}%
\definecolor{mycolor4}{rgb}{0.12941,0.12941,0.12941}%
\begin{tikzpicture}

\begin{axis}[%
  width=.4\linewidth,
  height=4cm,
  at={(0in,0.555in)},
  scale only axis,
  x dir=reverse,
  xmode=log,
  ymode=log,
  xmin=0.0196078431372549,
  xmax=0.166666666666667,
  ymin=1,
  ymax=10000,
  xminorticks=true,
  yminorticks=true,
  xlabel={$h$},
  ylabel={condition number},
  axis background/.style={fill=white},
  xmajorgrids, xminorgrids,
  ymajorgrids, yminorgrids,
  legend style={at={(0.03,1.1)}, anchor=north west, legend columns=3, legend cell align=left, align=left},
  legend image post style={scale=0.7},
  legend style={font=\small},
xtick={0.166666666666667,0.0909090909090909,0.0476190476190476,0.032258064516129,0.0196078431372549},
  xticklabels={0.17,0.09,0.05,0.03,0.02},
]
\addplot [color=mycolor1, line width=3.0pt, mark size=3.0pt, mark=o, mark options={solid, mycolor1}]
  table[row sep=crcr]{%
0.166666666666667	20\\
0.0909090909090909	49.0229566171763\\
0.0476190476190476	197.262230250342\\
0.032258064516129	443.810851776493\\
0.024390243902439	788.928571369973\\
0.0196078431372549	1233.04681935294\\
};
\addlegendentry{$A=H+S$}

\addplot [color=mycolor2, dashed, line width=3.0pt, mark size=1.8pt, mark=square, mark options={solid, mycolor2}]
  table[row sep=crcr]{%
0.166666666666667	20\\
0.0909090909090909	50\\
0.0476190476190476	199.999999999999\\
0.032258064516129	449.999999999999\\
0.024390243902439	799.999999999989\\
0.0196078431372549	1250\\
};
\addlegendentry{$H$}

\addplot [color=mycolor3, dashed, line width=3.0pt, mark size=3.0pt, mark=x, mark options={solid, mycolor3}]
  table[row sep=crcr]{%
0.166666666666667	2.6180339887499\\
0.0909090909090909	3.07768353717525\\
0.0476190476190476	6.31375151467505\\
0.032258064516129	9.51436445422258\\
0.024390243902439	12.7062047361747\\
0.0196078431372549	15.8945448438653\\
};
\addlegendentry{$S$}

\addplot [color=black, line width=3.0pt]
  table[row sep=crcr]{%
0.166666666666667	6\\
0.0909090909090909	11\\
0.0476190476190476	21\\
0.032258064516129	31\\
0.024390243902439	41\\
0.0196078431372549	51\\
};
\addlegendentry{$h^{-1}$}

\addplot [color=black, dashed, line width=3.0pt]
  table[row sep=crcr]{%
0.166666666666667	36\\
0.0909090909090909	121\\
0.0476190476190476	441\\
0.032258064516129	961\\
0.024390243902439	1681\\
0.0196078431372549	2601\\
};
\addlegendentry{$h^{-2}$}

\end{axis}

\begin{axis}[%
width=.4\linewidth,
height=4cm,
at={(.4\linewidth+1cm,0.555in)},
scale only axis,
x dir=reverse,
xmode=log,
xmin=0.0196078431372549,
xmax=0.166666666666667,
xminorticks=true,
xlabel style={font=\color{mycolor4}},
xlabel={$h$},
ymode=log,
ymin=1,
ymax=100,
yminorticks=true,
ylabel style={font=\color{mycolor4}},
axis background/.style={fill=white},
xmajorgrids,
xminorgrids,
ymajorgrids,
yminorgrids,
legend style={at={(0.03,1.1)}, anchor=north west,legend columns=3, legend cell align=left, align=left},
legend image post style={scale=0.7},
legend style={font=\small},
xtick={0.166666666666667,0.0909090909090909,0.0476190476190476,0.032258064516129,0.0196078431372549},
xticklabels={0.17,0.09,0.05,0.03,0.02},
]
\addplot [color=mycolor1, line width=3.0pt, mark size=1.7pt, mark=triangle, mark options={solid, mycolor1}]
  table[row sep=crcr]{%
0.166666666666667	1.69633928746495\\
0.0909090909090909	2.00788446887774\\
0.0476190476190476	2.18307013748588\\
0.032258064516129	2.24466188531328\\
0.024390243902439	2.27608709516835\\
0.0196078431372549	2.29514763766508\\
};
\addlegendentry{$H^{-1}A$}

\addplot [color=mycolor2, dashed, line width=3.0pt, mark size=4.9pt, mark=diamond, mark options={solid, mycolor2}]
  table[row sep=crcr]{%
0.166666666666667	1.69633928746495\\
0.0909090909090909	2.00788446887774\\
0.0476190476190476	2.18307013748588\\
0.032258064516129	2.24466188531329\\
0.024390243902439	2.27608709516835\\
0.0196078431372549	2.29514763766509\\
};
\addlegendentry{ichol$(H)^{-1}A$}

\addplot [color=black, line width=3.0pt]
  table[row sep=crcr]{%
0.166666666666667	6\\
0.0909090909090909	11\\
0.0476190476190476	21\\
0.032258064516129	31\\
0.024390243902439	41\\
0.0196078431372549	51\\
};
\addlegendentry{$h^{-1}$}

\end{axis}
\end{tikzpicture}
    \label{fig:cond_adv_diff}
\end{figure}

We now present a numerical illustration of Proposition~\ref{prop:adv_diff} and compute the condition numbers of a Galerkin discretization using piecewise affine linear finite elements. We set $\Omega = [0,1]\times [0,5] \times [0,1]$ with diffusivity $\nu = 0.001$, advection field $\mathbf {b}=(
    0.5,0,0
)^\top$ and $\mathbf{c}\equiv 0$. On the left-hand side of Figure~\ref{fig:cond_adv_diff}, we depict the condition numbers of the Galerkin projection of $\mathcal{A}$ denoted by $A$, of its symmetric part $\mathcal{H}$ denoted by $H$ and its skew-symmetric part $\mathcal{S}$ denoted by $S$. We observe in the left plot that the condition number of $A$ is dominated by that of the stiffness matrix of the diffusion term in the symmetric part behaving like $h^{-2}$ (due to the presence of two derivatives), while the condition number of the skew-symmetric part corresponding to the advection scales like $h^{-1}$, as it only involves one derivative. In the right plot, we observe that the condition number of the preconditioned matrix $H^{-1}A= I+H^{-1}S$ is bounded uniformly in the mesh width. This is due to the fact that, as proven in Proposition~\ref{prop:adv_diff}, $\calH^{-1}\calS$ is a bounded operator, hence has bounded spectrum, such that $I+\calH^{-1}\calS$ has a bounded spectrum as well. Correspondingly, the condition number for Galerkin projections is uniformly bounded.  The same can be observed  when approximating $H^{-1}$ with an incomplete Cholesky factorization using a fill-in tolerance of $10^{-2}$, see e.g. \cite{ScoT23}.

\subsection{Stokes equation}\label{subsec:stokes}
As a second example, we consider the stationary Stokes equation given by
\begin{equation}\label{eq:stokes}
\begin{split}
\begin{aligned}
    -\ddiv(\nu \nabla v) + \nabla P &= f&\qquad &\mathrm{on} \quad \Omega,\\
    \ddiv v &= 0 & \qquad &\mathrm{on} \quad \Omega, \\
  v &= 0& \qquad &\mathrm{on} \quad \partial\Omega 
\end{aligned}
\end{split}
\end{equation}
for $\Omega \subset \R^d$, $d\in \{2,3\}$. Here, $v:\Omega \to \R^d$ is the vector field modeling the velocity of the fluid and the scalar field $P:\Omega \to \R$ corresponds to the pressure. The differential operators $\nabla: \R^d \to \R^{d\times d},\ \operatorname{div}:\R^{d\times d}\to \R^d$ are to be understood row-wise and the kinematic viscosity $\nu: \R^d \to \R^{d\times d}$ is again a pointwise symmetric and uniformly positive matrix-valued mapping.

To include pressure-regularized formulations as commonly used in finite element methods, we introduce $s_1,s_2\in \R$ and consider the governing operator
\begin{align}\label{eq:stokes_stab}
    \mathcal{A} = 	\begin{bmatrix}
		-\ddiv(\nu \nabla v)  & \nabla \\
		\operatorname{div} & s_1 - s_2\Delta 
	\end{bmatrix}
\end{align}
which we may consider as an operator%
\begin{align*}
    \mathcal{A}&:L^2(\Omega;\R^d) \times L^2(\Omega) \supset D \to L^2(\Omega;\R^d) \times L^2(\Omega),\\
   D &= \{x\in H^1_0(\Omega;\R^d)\,|\, \nu \nabla x\in H(\ddiv,\Omega;\R^d)\} \times \{x \in H^1(\Omega)\,|\, s_2\Delta x\in L^2(\Omega)\}.
\end{align*}
For $s_2\neq 0$, we have $\{x \in H^1(\Omega)\,|\, s_2\Delta x\in L^2(\Omega)\} \cong H^2(\Omega)$ (here $ \cong$ denotes that the spaces are isomorphic) and for $s_2=0$ we have $\{x \in H^1(\Omega)\,|\, s_2\Delta x\in L^2(\Omega)\} = H^1(\Omega)$.
The operators then have the following structure.

\begin{proposition}%
\label{prop:strucstokes}
 The operator $\mathcal A$ in \eqref{eq:stokes_stab} may be decomposed into
    \begin{align*}
        \mathcal{H} = \begin{bmatrix}
         -\ddiv(\nu \nabla \cdot) & 0 \\
         0 & s_1 - s_2 \Delta 
        \end{bmatrix}, \qquad \mathcal{S} = \begin{bmatrix}
            0 & \nabla \\
            \ddiv & 0
        \end{bmatrix}
    \end{align*}
where $\calH:L^2(\Omega;\R^d)\times L^2(\Omega) \supset D \to L^2(\Omega;\R^d)\times L^2(\Omega)$ %
is symmetric and $\calS:L^2(\Omega;\R^d) \times L^2(\Omega)\supset D \to L^2(\Omega;\R^d)\times L^2(\Omega)$ %
is skew-symmetric. If $s_2>0$, then $\calH$ is self-adjoint.
\end{proposition}

\begin{proof}
Let us start with the top left block of $\calH$. The self-adjointness follows from the same arguments as in the proof of Lemma~\ref{lem:struc_advdiff}. For $s_2 \neq 0$, the same applies to the bottom right block. For $s_2=0$, we have the multiplication operator $s_1 I$ acting on $L^2(\Omega)$ which is clearly self-adjoint. %

 The skew-symmetry of $\calS$ follows, since for all $(v,P),(w,Q)\in D\subset H^1_0(\Omega;\R^d)\times H^1(\Omega)$ %
\begin{align}\label{eq:divgrad}
\begin{split}
        \langle \mathcal{S}(v,P),(w,Q)\rangle_{L^2(\Omega;\R^d)\times L^2(\Omega)} &=  \langle \nabla P, w\rangle_{L^2(\Omega;\R^d)} + \langle \operatorname{div} v, Q\rangle_{L^2(\Omega)} \\ 
    &= -\langle P,\ddiv w\rangle_{L^2(\Omega)} - \langle v,\nabla Q\rangle_{L^2(\Omega;\R^d)} %
\end{split}
\end{align}
due to the homogeneous Dirichlet boundary conditions.
\end{proof}

We briefly comment on the functional analytic framework of Proposition~\ref{prop:strucstokes}. The considered spaces only allow for skew-symmetry of $\calS$ which is in particular not closed on $D$. Loosely speaking, this results from domains which are not maximally chosen. Correspondingly, choosing suitable maximal domains, one may obtain also skew-adjointness. However, we stress that in usual finite element implementations \cite{braess1997efficient}, one would assemble a Galerkin discretization of $\calA$ on a finite element subspace $V\subset \operatorname{dom}(\calA)$ and then extract symmetric and skew-symmetric parts on the matrix level.  In this regard, we also note that one may consider a weak formulation in the sense of Remark~\ref{rem:weak} leading to bilinear forms. This will be illustrated  in the next subsection by means of the Oseen equation which is structurally similar to the Stokes equation.

We have the following result for the preconditioned operator.
\begin{proposition}
\label{prop:stokes}
Consider the operator $\mathcal A$ in \eqref{eq:stokes_stab}.
    Let $s_1\in \mathbb{R}\setminus \{0\}$ and set $X:= H(\ddiv;\Omega)\times H^1(\Omega)$. If $s_2 = 0$, then $\calH^{-1}\calS: X\to X$ is unbounded. If $s_2 > 0$, then $\calH^{-1}\calS: X\to X$ is bounded. In particular, $I + \calH^{-1}\calS\in L(X,X)$.
\end{proposition}
\begin{proof}
    Denote by $\calH_{11}$ and $\calH_{22}$ the diagonal blocks of $\calH$ endowed with the canonical domain inherited from the cross-product structure of $D = D_1\times D_2$. Then, $\calH_{11}$ is self-adjoint by a similar argument as in the proof of Proposition~\ref{prop:adv_diff} applied to vector-valued functions. Moreover, $\calH_{11}^{-1} \in L(L^2(\Omega;\R^d),D_1)$ which implies the boundedness of $\calH_{11}^{-1} \nabla \in L(H^1(\Omega),H^1(\Omega;\R^d))$ due to $\nabla \in L(H^1(\Omega),L^2(\Omega;\R^d)$ and $D_1\hookrightarrow H^1(\Omega;\R^d)$ (where $\hookrightarrow$ denotes the canonical embedding).

    Regarding the second block $\calH_{22}$, we distinguish the two cases. First, let $s_2=0$ such that $\calH_{22} = s_1 I$ and hence
    \begin{align*}
        \calH^{-1}\calS = \begin{bmatrix}
            0 & \calH_{11}^{-1}\nabla \\
            \tfrac{1}{s_1}\ddiv & 0.
        \end{bmatrix}
    \end{align*}
    The bottom left block is clearly unbounded as an operator from $H(\ddiv;\Omega)$ to $H^1(\Omega)$, as the weak divergence applied to an element from $H(\ddiv;\Omega)$ is (by definition of $H(\ddiv;\Omega)$) only square integrable and not necessarily weakly differentiable.

    In case $s_1 > 0$, $\calH_{22}$ is closed such that $\calH_{22} = s_1 - s_2\Delta$ is boundedly invertible (as a mapping from $L^2(\Omega)$ to $\dom(\calH_{22}) = H^2(\Omega)\hookrightarrow H^1(\Omega)$)  due to the Lax-Milgram theorem and a similar argument as in the proof of Proposition~\ref{prop:adv_diff}. In particular, $\calH_{22}^{-1}\ddiv\in L(H(\ddiv;\Omega),H^1(\Omega))$. Hence $\calH^{-1}\calS$ is bounded as all blocks are bounded.
\end{proof}

\begin{figure}[htb]
	\centering
\definecolor{mycolor1}{rgb}{0.06600,0.44300,0.74500}%
\definecolor{mycolor2}{rgb}{0.86600,0.32900,0.00000}%
\definecolor{mycolor3}{rgb}{0.92900,0.69400,0.12500}%
\definecolor{mycolor4}{rgb}{0.52100,0.08600,0.81900}%
\definecolor{mycolor5}{rgb}{0.23100,0.66600,0.19600}%
\definecolor{mycolor6}{rgb}{0.12941,0.12941,0.12941}%

\begin{tikzpicture}

\begin{axis}[%
width=.4\linewidth,
height=4cm,
at={(0,0.552in)},
scale only axis,
x dir=reverse,
xmode=log,
xmin=0.0196078431372549,
xmax=0.166666666666667,
xminorticks=true,
xlabel={$h$},
ymode=log,
ymin=360,
ymax=1000000,
ytick={   1000,   10000,  100000, 1000000},
yminorticks=true,
ylabel={condition number},
axis background/.style={fill=white},
xmajorgrids,
xminorgrids,
ymajorgrids,
yminorgrids,
legend style={at={(0.03,1.2)}, anchor=north west,legend columns=2, legend cell align=left, align=left},
legend image post style={scale=0.7},
legend style={font=\small},
xtick={0.166666666666667,0.0909090909090909,0.0476190476190476,0.032258064516129,0.0196078431372549},
xticklabels={0.17,0.09,0.05,0.03,0.02},
]
\addplot [color=mycolor1, line width=3.0pt, mark size=3.0pt, mark=o, mark options={solid, mycolor1}]
  table[row sep=crcr]{%
0.166666666666667	6382.49617388486\\
0.0909090909090909	21218.166822947\\
0.0476190476190476	76817.9158276425\\
0.032258064516129	166976.204338395\\
0.024390243902439	291685.660744887\\
0.0196078431372549	450942.082686409\\
};
\addlegendentry{$A$}

\addplot [color=mycolor2, dashed, line width=3.0pt, mark size=1.8pt, mark=square, mark options={solid, mycolor2}]
  table[row sep=crcr]{%
0.166666666666667	7311.68523111275\\
0.0909090909090909	24574.8949172392\\
0.0476190476190476	89566.3523743717\\
0.032258064516129	195177.470820343\\
0.024390243902439	341408.250207072\\
0.0196078431372549	528258.690534559\\
};
\addlegendentry{$H$}

\addplot [color=mycolor3, line width=3.0pt, mark size=1.7pt, mark=triangle, mark options={solid, mycolor3}]
  table[row sep=crcr]{%
0.166666666666667	828.426680930923\\
0.0909090909090909	2674.39407021129\\
0.0476190476190476	9489.48511805872\\
0.032258064516129	20455.5704049597\\
0.024390243902439	35561.6658969284\\
0.0196078431372549	54804.1367375154\\
};
\addlegendentry{$H^{-1}A$}

\addplot [color=mycolor4, dashed, line width=3.0pt, mark size=4.9pt, mark=diamond, mark options={solid, mycolor4}]
  table[row sep=crcr]{%
0.166666666666667	827.78172912857\\
0.0909090909090909	2672.09276333941\\
0.0476190476190476	9481.25215785972\\
0.032258064516129	23395.6139059211\\
0.024390243902439	54699.2436948438\\
0.0196078431372549	105924.325657722\\
};
\addlegendentry{ichol$(H,10^{-2})^{-1}A$}

\addplot [color=black, dashed, line width=3.0pt]
  table[row sep=crcr]{%
0.166666666666667    360\\
0.0909090909090909   1210\\
0.0476190476190476   4410\\
0.032258064516129    9610\\
0.024390243902439    16810\\
0.0196078431372549   26010\\
};
\addlegendentry{$h^{-2}$}
\end{axis}

\begin{axis}[%
width=.4\linewidth,
height=4cm,
at={(.4\linewidth+1cm,0.555in)},
scale only axis,
x dir=reverse,
xmode=log,
xmin=0.0196078431372549,
xmax=0.166666666666667,
xminorticks=true,
xlabel style={font=\color{mycolor6}},
xlabel={$h$},
ymode=log,
ymin=1,
ymax=100000,
ytick={   10, 100,1000,   10000, 1000000,10^5},
yminorticks=true,
ylabel style={font=\color{mycolor6}},
axis background/.style={fill=white},
xmajorgrids,
xminorgrids,
ymajorgrids,
yminorgrids,
legend style={at={(0.97,0.5)}, anchor=east, legend columns=6, legend cell align=left, align=left},
legend style={at={(0.03,1.2)}, anchor=north west,legend columns=2, legend cell align=left, align=left},
legend image post style={scale=0.7},
legend style={font=\small},
xtick={0.166666666666667,0.0909090909090909,0.0476190476190476,0.032258064516129,0.0196078431372549},
xticklabels={0.17,0.09,0.05,0.03,0.02},
]
\addplot [color=mycolor1, line width=3.0pt, mark size=3.0pt, mark=o, mark options={solid, mycolor1}]
  table[row sep=crcr]{%
0.166666666666667	433.474693041144\\
0.0909090909090909	1228.68891993496\\
0.0476190476190476	4054.29071399265\\
0.032258064516129	8520.09667780102\\
0.024390243902439	14626.234206432\\
0.0196078431372549	22372.745794436\\
};
\addlegendentry{$A$}

\addplot [color=mycolor2, dashed, line width=3.0pt, mark size=1.8pt, mark=square, mark options={solid, mycolor2}]
  table[row sep=crcr]{%
0.166666666666667	404.349379884933\\
0.0909090909090909	1171.81748870093\\
0.0476190476190476	3905.94557554142\\
0.032258064516129	8239.68648828502\\
0.024390243902439	14173.2566209601\\
0.0196078431372549	21706.7298278856\\
};
\addlegendentry{$H$}

\addplot [color=mycolor3, line width=3.0pt, mark size=1.7pt, mark=triangle, mark options={solid, mycolor3}]
  table[row sep=crcr]{%
0.166666666666667	2.48220554372606\\
0.0909090909090909	2.3902873248791\\
0.0476190476190476	2.29773651438766\\
0.032258064516129	2.28691888521347\\
0.024390243902439	2.26821822891063\\
0.0196078431372549	2.25548583896427\\
};
\addlegendentry{$H^{-1}A$}

\addplot [color=mycolor4, dashed, line width=3.0pt, mark size=4.9pt, mark=diamond, mark options={solid, mycolor4}]
  table[row sep=crcr]{%
0.166666666666667	19.9026879710886\\
0.0909090909090909	110.933822743966\\
0.0476190476190476	343.502349156802\\
0.032258064516129	718.435270997818\\
0.024390243902439	1230.45493982978\\
0.0196078431372549	1878.46505312403\\
};
\addlegendentry{ichol$(H,10^{-2})^{-1}A$}
\addplot [color=black, dashed, line width=3.0pt]
  table[row sep=crcr]{%
0.166666666666667	3.6\\
0.0909090909090909	12.1\\
0.0476190476190476	44.1\\
0.032258064516129	96.1\\
0.024390243902439	168.1\\
0.0196078431372549	260.1\\
};
\addlegendentry{$h^{-2}$}
\addplot [color=mycolor5, dashed, line width=3.0pt, mark size=4.9pt, mark=diamond, mark options={solid, mycolor5}]
  table[row sep=crcr]{%
0.166666666666667	2.48220554372606\\
0.0909090909090909	2.15765384613053\\
0.0476190476190476	4.79113042804431\\
0.032258064516129	10.8594758918581\\
0.024390243902439	20.2316573556474\\
0.0196078431372549	34.8111584083132\\
};
\addlegendentry{ichol$(H,10^{-4})^{-1}A$}

\end{axis}
\end{tikzpicture}%
 	\caption{Condition numbers for Stokes equation with $L^2$ pressure regularization $(s_1,s_2) = (1,0)$ (left) and $H^1$ pressure regularization  $(s_1,s_2) = (1,1)$ (right).}
	\label{fig:stokes_cond}
\end{figure}
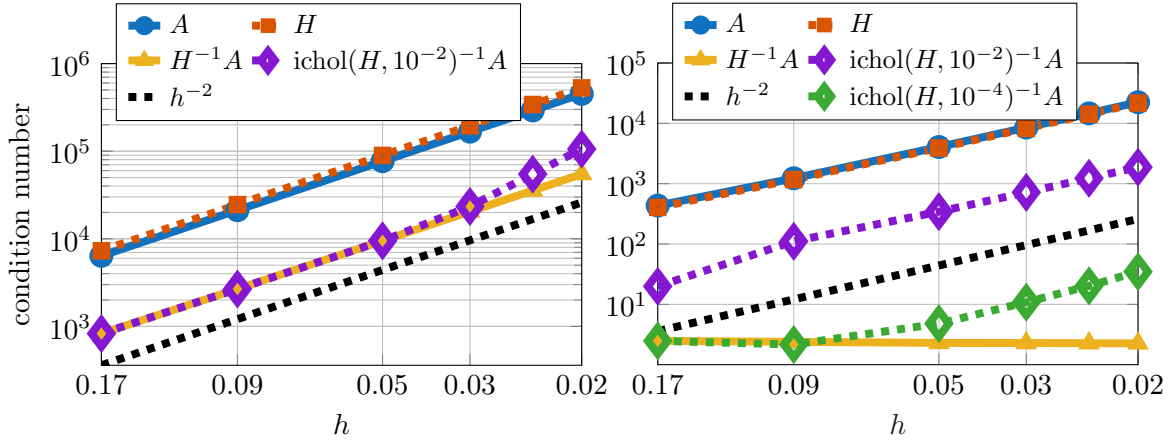

The numerical illustration of Proposition~\ref{prop:stokes} is presented in Figure~\ref{fig:stokes_cond}, where we chose $\Omega = [0,1]^2$, $\nu \equiv 1$ and discretized the system with P2 vector-valued finite elements for the velocity field $v:\Omega \to \R^2$ and P1 finite elements for the pressure $P: \Omega \to \R$. We note that for this choice of finite elements, one usually does not employ a stabilization. Here, we do so regardless, to obtain invertibility of $\calH$. The unstabilized case will be considered in the next subsection. 

In the left plot of Figure~\ref{fig:stokes_cond}, we observe an increasing condition number if we only employ an $L^2$-pressure stabilization, corresponding to the choice $s_2=0$ in Proposition~\ref{prop:stokes}. This is due to the fact that $\calH^{-1}\calS$ contains a (first-order) differential operator in the bottom left block, leading to a behavior of the condition number like $\mathcal{O}(h^{-1})$. In the right plot of Figure~\ref{fig:stokes_cond}, we see that a $H^1$-regularization ($s_2>0$) leads to a uniformly bounded  (in the mesh size) condition number, which reflects the boundedness of $\calH^{-1}\calS$ as proven in Proposition~\ref{prop:stokes}. Furthermore, the inversion of the block-wise elliptic operator $H$ may be approximated with an incomplete Cholesky factorization, where we however observe an increase of the condition number for small mesh sizes (depending on the fill-in tolerance).  Note that this example was also considered in \cite[Examples 7.1 and 7.2]{mardal2011preconditioning}, where the authors observed a similar behavior of the condition number and discussed suitable inf-sup conditions for finite element discretizations.

\subsection{Oseen equation}\label{subsec:oseen}
As third example we consider the Oseen equation given by
\begin{equation}\label{eq:oseen}
\begin{split}
\begin{aligned}
     -\ddiv \sigma(v,P) + \mathbf{b}^\top \nabla v &= f &\qquad &\mathrm{on} \quad \Omega,\\
    \ddiv v &= 0   & \qquad &\mathrm{on} \quad \Omega, \\
  v &= 0  & \qquad &\mathrm{on} \quad \partial\Omega ,
\end{aligned}
\end{split}
\end{equation}
which arises when linearizing the Navier-Stokes equation at a stationary profile. Again, we choose $\Omega \subset \R^d$, $d\in \{2,3\}$ and with the dynamic viscosity $\mu > 0$, we denote the symmetric stress tensor by
\[
\sigma(v,P) = \mu( \nabla v + \nabla v^\top)-PI_d.
\]
Writing \eqref{eq:oseen} in a block-wise fashion, we consider the operator
\begin{equation*}
    \mathcal{A}=\begin{bmatrix} \calA_{11} & \calA_{12}\\ \calA_{21} & 0\end{bmatrix}:= \begin{bmatrix}
        -\ddiv\nu (\nabla \cdot  + \nabla^\top\cdot ) + \mathbf{b}^\top \nabla \cdot & \nabla \\
        \ddiv & 0
    \end{bmatrix}
\end{equation*}
as a mapping
\begin{align*}
   \calA &: H^{-1}(\Omega;\R^d) \times H^{-1}(\Omega) \supset \dom(\calA) \to H^{-1}(\Omega;\R^d) \times H^{-1}(\Omega)\\
\dom(\calA) &= H^1_0(\Omega;\R^d) \times L^2(\Omega).
\end{align*}
Note that here (and in contrast to the advection-diffusion-reaction equation in Subsection~\ref{subsec:advec}), we interpret the second-order differential  operator in divergence form in the top left block in weak form, i.e., for $v,w\in H^1(\Omega;\R^d)$,
\begin{equation}\label{eq:defA11}
    \langle \mathcal{A}_{11}v,w\rangle_{H^{-1}(\Omega;\R^d), H^1(\Omega;\R^d)} = \langle \nu ( \nabla v + \nabla^\top v), \nabla w\rangle_{L^2(\Omega;\R^{d\times d})} + \langle \mathbf{b}^\top \nabla v, w\rangle_{L^2(\Omega;\R^d)}.
\end{equation}
Similarly, $\calA_{12}$ is to be understood as operator
\begin{equation}\label{eq:defA12}
  \langle \calA_{12} P, w\rangle_{H^{-1}(\Omega;\R^d), H^{1}(\Omega;\R^d)} = -\langle P,\ddiv w\rangle_{L^2(\Omega)} \  \mbox{\rm for all } P\in H^1(\Omega), w\in H^1_0(\Omega).
\end{equation}
In contrast to the Stokes equation considered in Subsection~\ref{subsec:stokes}, we do not employ a pressure regularization here. Hence, the bottom right block of $\calA$ is zero, such that the symmetric part is not invertible. To deal with this case, we consider the Schur complement approach as suggested in \cite{GuduLies22}, i.e., we consider the Schur complement of $\calA_{11}$ in $\calA$ given by
\begin{equation}\label{eq:schuroseen}
    \calW = -\calA_{21} \calA_{11}^{-1}\calA_{12}. %
\end{equation}
and obtain the following structure and boundedness result.

\begin{proposition}
\label{prop:oseen}
The Schur complement $\calW: L^2(\Omega) \supset H^1(\Omega) \to L^2(\Omega)$ is bounded, that is, $\calW\in L(L^2(\Omega), L^2(\Omega))$ and accretive.
\end{proposition}
\begin{proof}
First we show that $\calA_{11}$ as defined in \eqref{eq:defA11} induces a coercive and bounded bilinear form in $H^1_0(\Omega;\R^d)$. The boundedness is clear by definition.
To show the coercivity, note that for all $v \in H_0^1(\Omega;\R^d)$ we have
\begin{equation}\label{eq:advanishes}
    \langle \mathbf {b}^\top \nabla v,v\rangle_{L^2(\Omega;\R^d)} = - \langle v,\mathbf {b}^\top \nabla v\rangle_{L^2(\Omega;\R^d)}
\end{equation}
as $\mathbf b$ is divergence free by definition
and due to the no slip boundary condition $v=0$. Thus, $\langle \mathbf b^\top \nabla v,v\rangle_{L^2(\Omega;\R^d)} = 0$. Moreover, by Korn's inequality \cite[Theorem 6.3.4]{Ciar21}, there exists a constant $c>0$ such that
\begin{equation*}
    \langle \nu (\nabla v + \nabla^\top v), \nabla v\rangle_{L^2(\Omega;\R^{d\times d})} \geq c\|v\|^2_{H^1(\Omega;\R^d)}.
\end{equation*}
Together with \eqref{eq:advanishes}, this implies that the bilinear form is coercive, i.e.
\begin{equation*}
   \langle \calA_{11} v,v\rangle_{H^{-1}(\Omega;\R^d),H^1(\Omega;\R^d)}\geq c\|v\|_{H^1(\Omega;\R^d)}^2 \qquad \mbox{\rm for all } v\in H^1_0(\Omega;\R^{d}).
\end{equation*}
Thus, 
$\calA_{11}$ is invertible and $\calA_{11}^{-1}:H^{-1}(\Omega;\R^d) \to H^1_0(\Omega;\R^d)$.
Hence, in view of \eqref{eq:defA12},
\begin{equation*}
    \calA_{11}^{-1}\calA_{12} \in L(L^2(\Omega), H^1(\Omega;\R^d)).
\end{equation*}
such that, applying the divergence operator $\calA_{21} = \ddiv \in L(H^1(\Omega;\R^d), L^2(\Omega))$, we obtain boundedness of the Schur complement, that is $\calW = \calA_{21} \calA_{11}^{-1}\calA_{12} \in L(L^2(\Omega),L^2(\Omega))$.

To see that $\calW$ is accretive, we show that $\calA_{21} = \calA_{12}^*$ (where we consider the former as an operator $\calA_{21}\in L(H^1(\Omega;\R^d),L^2(\Omega))$ and $\calA_{12}\in L(L^2(\Omega)),H^{-1}(\Omega;\R^d))$ and that $\calA_{11}^{-1}$ is accretive. This follows due to
\begin{equation*}
    \langle \calA_{12} P, w\rangle_{H^{-1}(\Omega;\R^d), H^{1}(\Omega;\R^d)} = -\langle P,\ddiv w\rangle_{L^2(\Omega)} = -\langle P,\calA_{21}w\rangle_{L^2(\Omega)}
\end{equation*}
by the no-slip boundary condition $w=0$ on $\partial \Omega$. For the accretivity of $\calA_{11}^{-1}$, let $v\in H^{-1}(\Omega;\R^d)$. Then, by the invertiblity of $\calA_{11}: H_0^1(\Omega;\R^d)\to H^{-1}(\Omega;\R^d)$ there is a unique $w\in H^1_0(\Omega;\R^d)$ such that $\calA_{11} w = v$. Hence,
\begin{align*}
    \langle \calA_{11}^{-1} v,v\rangle_{H^1(\Omega;\R^d),H^{-1}(\Omega;\R^d)} &= \langle \calA_{11}^{-1} \calA_{11}w,\calA_{11}w\rangle_{H^1(\Omega;\R^d),H^{-1}(\Omega;\R^d)} \\
    &= \langle w,\calA_{11}w\rangle_{H^1(\Omega;\R^d),H^{-1}(\Omega;\R^d)} \\
    &\geq c\|w\|_{H^1(\Omega;\R^d)}^2 = c\|\calA_{11}^{-1}v\|_{H^1(\Omega;\R^d)}^2 \geq \tilde c \|v\|_{H^{-1}(\Omega;\R^d)}^2,
\end{align*}
where the second last inequality follows from accretivity of $\calA_{11}$ and the last inequality follows with $\tilde c = \frac{c}{\overline{c}}$ from the boundedness, i.e., that $\|\calA_{11}\|_{L(H^1(\Omega;\R^d),H^{-1}(\Omega;\R^d))}\leq \overline{c}$.
Consequently,
\begin{align*}
    \langle \calW v,v\rangle_{L^2(\Omega)} &= -\langle \calA_{21}\calA_{11}^{-1}\calA_{12} v,v\rangle_{L^2(\Omega)} \\
    &= -\langle \calA_{11}^{-1}\calA_{12} v, \calA_{21}^* v\rangle_{H^1(\Omega;\R^d),H^{-1}(\Omega;\R^d)} \\
    &= \langle \calA_{11}^{-1}\calA_{12} v, \calA_{12} v\rangle_{H^1(\Omega;\R^d),H^{-1}(\Omega;\R^d)} \geq \tilde c\|\calA_{12}v\|^2_{H^{-1}(\Omega;\R^d)}\geq 0,
\end{align*}
which shows the accretivity of $\calW$.
\end{proof}

We briefly comment on the last part of the proof of Proposition~\ref{prop:oseen}. When considering $L^2(\Omega)$ functions with zero mean, we may deduce strict accretivity by means of an extension of the Poincaré inequality $\|\calA_{12}v\|_{H^{-1}(\Omega;\R^d)} = \|\nabla v\|_{H^{-1}(\Omega;\R^d)} \geq C \|v\|_{L^2(\Omega;\R^d)}$ for some constant $C>0$ and all $v\in L^2(\Omega;\R^d)$.

\begin{figure}[htb]
	\centering
\definecolor{mycolor1}{rgb}{0.06600,0.44300,0.74500}%
\definecolor{mycolor2}{rgb}{0.86600,0.32900,0.00000}%
\definecolor{mycolor3}{rgb}{0.92900,0.69400,0.12500}%
\definecolor{mycolor4}{rgb}{0.12941,0.12941,0.12941}%
\begin{tikzpicture}

\begin{axis}[%
  width=.4\linewidth,
  height=4cm,
  at={(0in,0.555in)},
  scale only axis,
  x dir=reverse,
  xmode=log,
  ymode=log,
  xmin= 0.0384615384615385,
  xmax=0.166666666666667,
  ymin=1,
  ymax=10000,
  xminorticks=true,
  yminorticks=true,
  xlabel={$h$},
  ylabel={condition number},
  axis background/.style={fill=white},
  xmajorgrids, xminorgrids,
  ymajorgrids, yminorgrids,
  legend style={at={(0.03,1.1)}, anchor=north west, legend columns=3, legend cell align=left, align=left},
  legend image post style={scale=0.7},
  legend style={font=\small},
  xtick={0.166666666666667,0.0909090909090909,0.0476190476190476,0.032258064516129},
  xticklabels={0.17,0.09,0.05,0.03},
]

\addplot [color=mycolor1, line width=3.0pt, mark size=3.0pt, mark=o, mark options={solid, mycolor1}]
  table[row sep=crcr]{%
0.166666666666667 114.296640838772\\
0.0909090909090909 380.448220633856\\
0.0666666666666667 705.264973422568\\
0.0476190476190476 1377.85750381368\\
0.0384615384615385 2112.19358376743\\
};
\addlegendentry{$A_{11}$}

\addplot [color=mycolor2, dashed, line width=3.0pt, mark size=1.8pt, mark=square, mark options={solid, mycolor2}]
  table[row sep=crcr]{%
0.166666666666667 115.391621977278\\
0.0909090909090909 379.283368611903\\
0.0666666666666667 706.068009708717\\
0.0476190476190476 1384.95120857922\\
0.0384615384615385 2126.27149012571\\
};
\addlegendentry{$H_{11}$}

\addplot [color=mycolor3, dashed, line width=3.0pt, mark size=3.0pt, mark=x, mark options={solid, mycolor3}]
  table[row sep=crcr]{%
0.166666666666667 1.9400411355703\\
0.0909090909090909 2.04142620707728\\
0.0666666666666667 2.08449509487298\\
0.0476190476190476 2.11743360799054\\
0.0384615384615385 2.13568237568791\\
};
\addlegendentry{$H_{11}^{-1} A_{11}$}

\addplot [color=black, line width=3.0pt]
  table[row sep=crcr]{%
0.166666666666667 6\\
0.0909090909090909 11\\
0.0666666666666667 15\\
0.0476190476190476 21\\
0.0384615384615385 26\\
};
\addlegendentry{$h^{-1}$}

\addplot [color=black, dashed, line width=3.0pt]
  table[row sep=crcr]{%
0.166666666666667 36\\
0.0909090909090909 121\\
0.0666666666666667 225\\
0.0476190476190476 441\\
0.0384615384615385 676\\
};
\addlegendentry{$h^{-2}$}

\end{axis}

\begin{axis}[%
  width=.4\linewidth,
  height=4cm,
  at={(.4\linewidth+1cm,0.555in)},
  scale only axis,
  x dir=reverse,
  xmode=log,
  xmin=0.0384615384615385,
  xmax=0.166666666666667,
  xminorticks=true,
  xlabel style={font=\color{mycolor4}},
  xlabel={$h$},
  ymode=log,
  ymin=1,
  ymax=1000,
  yminorticks=true,
  axis background/.style={fill=white},
  xmajorgrids,
  xminorgrids,
  ymajorgrids,
  yminorgrids,
  legend style={at={(0.03,1.1)}, anchor=north west, legend columns=3, legend cell align=left, align=left},
  legend image post style={scale=0.7},
  legend style={font=\small},
  xtick={0.166666666666667,0.0909090909090909,0.0476190476190476,0.032258064516129},
  xticklabels={0.17,0.09,0.05,0.03},
]

\addplot [color=mycolor1, line width=3.0pt, mark size=1.7pt, mark=triangle, mark options={solid, mycolor1}]
  table[row sep=crcr]{%
0.166666666666667 72.2482611527849\\
0.0909090909090909 87.1539023622475\\
0.0666666666666667 91.3410327440329\\
0.0476190476190476 94.1711306135971\\
0.0384615384615385 96.1334011876496\\
};
\addlegendentry{$W$}

\addplot [color=mycolor2, dashed, line width=3.0pt, mark size=4.9pt, mark=diamond, mark options={solid, mycolor2}]
  table[row sep=crcr]{%
0.166666666666667 1.01173353397136\\
0.0909090909090909 1.01202993930035\\
0.0666666666666667 1.01204826915966\\
0.0476190476190476 1.0120535065481\\
0.0384615384615385 1.01205453385887\\
};
\addlegendentry{$H_W^{-1}\,W$}

\addplot [color=mycolor3, dashed, line width=3.0pt, mark size=3.0pt, mark=x, mark options={solid, mycolor3}]
  table[row sep=crcr]{%
0.166666666666667 11.9328452960204\\
0.0909090909090909 11.7456472330874\\
0.0666666666666667 11.6853434809876\\
0.0476190476190476 11.6376463482373\\
0.0384615384615385 11.614571080857\\
};
\addlegendentry{$M_p^{-1}\,W$}

\addplot [color=black, line width=3.0pt]
  table[row sep=crcr]{%
0.166666666666667 6\\
0.0909090909090909 11\\
0.0666666666666667 15\\
0.0476190476190476 21\\
0.0384615384615385 26\\
};
\addlegendentry{$h^{-1}$}

\addplot [color=black, dashed, line width=3.0pt]
  table[row sep=crcr]{%
0.166666666666667 36\\
0.0909090909090909 121\\
0.0666666666666667 225\\
0.0476190476190476 441\\
0.0384615384615385 676\\
};
\addlegendentry{$h^{-2}$}

\end{axis}
\end{tikzpicture}%
 	\caption{Condition numbers for Oseen equation for the original system (left) and the Schur complement (right).}
	\label{fig:oseen_b}
\end{figure}
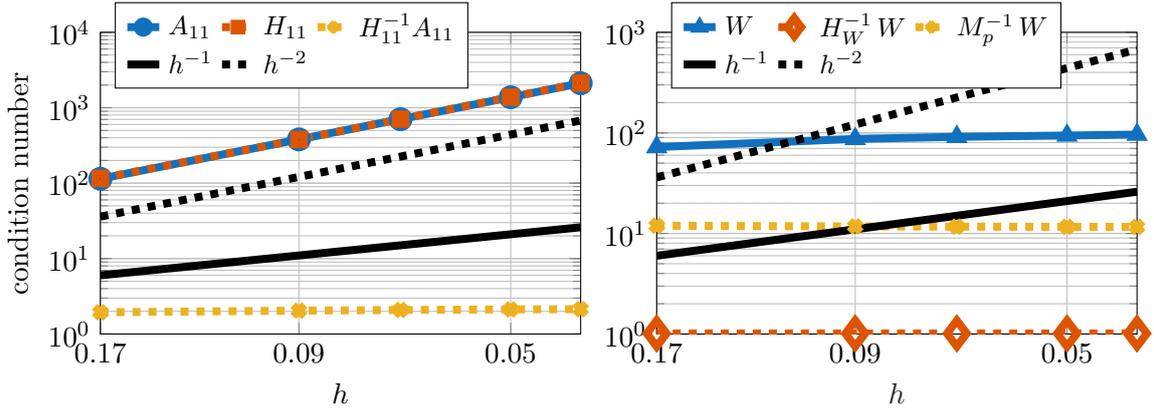

We depict the condition numbers associated with the Oseen equation in Figure~\ref{fig:oseen_b}. We choose the same parameters and domain as for the Stokes equation of Subsection~\ref{subsec:stokes} and the advection term $\mathbf b=[1,1]^\top$. On the left-hand side of Figure~\ref{fig:oseen_b}, and as to be expected, we observe the same behavior for the top left block as for the advection-diffusion equation, that is quadratic increase of the condition number for decreasing mesh width. On the right hand side of Figure~\ref{fig:oseen_b} we depict the condition number of the Schur complement $W$ (being a discretization of $\calW$ from \eqref{eq:schuroseen}) and its preconditioned versions. We observe that already the unpreconditioned Schur complement has a uniformly bounded condition number (due to its boundedness proven in Proposition~\ref{prop:oseen}, while a preconditioning with a pressure mass matrix $M_p$ decreases this condition number by one order of magnitude, and preconditioning with the symmetric part of the Schur complement $H_W$ decreases it by another order of magnitude. Note that, however, for the latter it is necessary to form the Schur complement, while a preconditioning with the mass matrix of the pressure can be done very efficiently, e.g.\ only using its diagonal part~\cite{wathen1987realistic}. %

\subsection{Wave equation with momentum damping}\label{subsec:wave}
As a next example, we consider a hyperbolic PDE given by a wave equation. %
The stationary system corresponding to the first order port-Hamiltonian formulation reads 
\begin{equation*}
\begin{aligned}
	-\ddiv q + \rho p &= f_1& \qquad &\mathrm{on} \quad \Omega,\\
	-\nabla p + \eta q &= f_2& \qquad &\mathrm{on} \quad \Omega,\\
    	n^\top q &= 0 &\qquad \ \, &\mathrm{on} \quad \partial\Omega.
\end{aligned}
\end{equation*}
Here, $p\colon  \Omega \to \R$ is the momentum, $q\colon \Omega \to \R^d$ is the vector-valued strain and $\rho,\eta \geq 0$ are scalar friction parameters.  The boundary condition models that there is no strain in normal direction. The governing operator of this system is given by
\begin{equation}\label{waveop}
\calA = 
	\begin{bmatrix}
		\rho & -\ddiv \\
		-\nabla & \eta
	\end{bmatrix}
\end{equation}
which we consider as a mapping
\begin{equation*}
\calA : L^2(\Omega)\times L^2(\Omega;\R^d) \supset H^1(\Omega) \times H(\ddiv;\Omega) \to L^2(\Omega)\times L^2(\Omega;\R^d).
\end{equation*}

\begin{proposition}\label{eq:wave_decomp}
    The operator $\calA$ in \eqref{waveop} may be decomposed as
    \begin{equation*}
  \calH+\calS=    \begin{bmatrix}
	    \rho & 0 \\
        0 &\eta
	\end{bmatrix}+ \begin{bmatrix}
		0 & -\ddiv \\
		-\nabla & 0
	\end{bmatrix}, 
    \end{equation*}
where $\calS:L^2(\Omega)\times L^2(\Omega;\R^d) \supset H^1(\Omega) \times H(\ddiv;\Omega) \to L^2(\Omega)\times L^2(\Omega;\R^d)$ is skew-adjoint and $\calH: L^2(\Omega)\times L^2(\Omega;\R^d)\to L^2(\Omega)\times L^2(\Omega;\R^d)$ is self-adjoint.
\end{proposition}
\begin{proof}
 The self-adjointness of $\calH$ as multiplication operator is obvious. The skew-adjointness of $\calS$ follows from the analogous computation to \eqref{eq:divgrad}.
\end{proof}

The next proposition follows straightforwardly, as the symmetric part is a bounded operator, and hence its inverse cannot compensate for the unboundedness of the derivative operators.
\begin{proposition}
\label{prop:wave} Consider the operator
$\calA$ in \eqref{waveop} and its decomposition from Proposition~\ref{eq:wave_decomp}.
    Let $X\in \{L^2(\Omega)\times L^2(\Omega;\R^d), H^1(\Omega) \times H(\ddiv;\Omega)  \}$. Then $\calH^{-1}\calS: X\to X$ is unbounded.
\end{proposition}

\begin{figure}[htb]
	\centering
\definecolor{mycolor1}{rgb}{0.06600,0.44300,0.74500}%
\definecolor{mycolor2}{rgb}{0.86600,0.32900,0.00000}%
\definecolor{mycolor3}{rgb}{0.12941,0.12941,0.12941}%
\definecolor{mycolor4}{rgb}{0.92900,0.69400,0.12500}%

\begin{tikzpicture}

\begin{axis}[%
  width=.4\linewidth,
  height=4cm,
  at={(0in,0.555in)},
  scale only axis,
  x dir=reverse,
  xmode=log,
  xmin=0.0196078431372549,
  xmax=0.0476190476190476,
  xminorticks=true,
  xlabel style={font=\color{mycolor3}},
  xlabel={$h$},
  ymode=log,
  ymin=1,
  ymax=10000,
  yminorticks=true,
  ylabel style={font=\color{mycolor3}},
  ylabel={condition number},
  axis background/.style={fill=white},
  xmajorgrids, xminorgrids,
  ymajorgrids, yminorgrids,
  legend style={at={(0.03,1.1)}, anchor=north west, legend columns=3, legend cell align=left, align=left},
  legend image post style={scale=0.7},
  legend style={font=\small},
  xtick={0.166666666666667,0.0909090909090909,0.0476190476190476,0.032258064516129,0.0196078431372549},
  xticklabels={0.17,0.09,0.05,0.03,0.02},
]

\addplot [color=mycolor1, line width=3.0pt, mark size=3.0pt, mark=o, mark options={solid, mycolor1}]
  table[row sep=crcr]{%
0.0476190476190476  80.5492027763024\\
0.032258064516129   114.35247313722\\
0.024390243902439   150.420667233851\\
0.0196078431372549  185.144976311951\\
};
\addlegendentry{$A=H+S$}

\addplot [color=mycolor2, dashed, line width=3.0pt, mark size=1.8pt, mark=square, mark options={solid, mycolor2}]
  table[row sep=crcr]{%
0.0476190476190476  4.73205080757226\\
0.032258064516129   4.73205080756888\\
0.024390243902439   4.73205080756888\\
0.0196078431372549  4.73205080756889\\
};
\addlegendentry{$H$}

\addplot [color=black, line width=3.0pt]
  table[row sep=crcr]{%
0.0476190476190476  21\\
0.032258064516129   31\\
0.024390243902439   41\\
0.0196078431372549  51\\
};
\addlegendentry{$h^{-1}$}

\addplot [color=black, dashed, line width=3.0pt]
  table[row sep=crcr]{%
0.0476190476190476  441\\
0.032258064516129   961\\
0.024390243902439   1681\\
0.0196078431372549  2601\\
};
\addlegendentry{$h^{-2}$}

\end{axis}

\begin{axis}[%
  width=.4\linewidth,
  height=4cm,
  at={(.4\linewidth+1cm,0.555in)},
  scale only axis,
  x dir=reverse,
  xmode=log,
  xmin=0.0196078431372549,
  xmax=0.0476190476190476,
  xminorticks=true,
  xlabel style={font=\color{mycolor3}},
  xlabel={$h$},
  ymode=log,
  ymin=21,
  ymax=309.422399718753,
  yminorticks=true,
  ylabel style={font=\color{mycolor3}},
  axis background/.style={fill=white},
  xmajorgrids, xminorgrids,
  ymajorgrids, yminorgrids,
  legend style={at={(0.03,1.1)}, anchor=north west, legend columns=3, legend cell align=left, align=left},
  legend image post style={scale=0.7},
  legend style={font=\small},
  xtick={0.166666666666667,0.0909090909090909,0.0476190476190476,0.032258064516129,0.0196078431372549},
  xticklabels={0.17,0.09,0.05,0.03,0.02},
]

\addplot [color=mycolor1, line width=3.0pt, mark size=1.7pt, mark=triangle, mark options={solid, mycolor1}]
  table[row sep=crcr]{%
0.0476190476190476  129.653567854746\\
0.032258064516129   189.579059928144\\
0.024390243902439   249.501420507495\\
0.0196078431372549  309.422399718748\\
};
\addlegendentry{$H^{-1}A$}

\addplot [color=mycolor2, dashed, line width=3.0pt, mark size=4.9pt, mark=diamond, mark options={solid, mycolor2}]
  table[row sep=crcr]{%
0.0476190476190476  129.653567854746\\
0.032258064516129   189.579059928145\\
0.024390243902439   249.501420507507\\
0.0196078431372549  309.422399718753\\
};
\addlegendentry{ichol$(H)^{-1}A$}

\addplot [color=black, line width=3.0pt]
  table[row sep=crcr]{%
0.0476190476190476  21\\
0.032258064516129   31\\
0.024390243902439   41\\
0.0196078431372549  51\\
};
\addlegendentry{$h^{-1}$}

\end{axis}
\end{tikzpicture}%
 	\caption{Condition numbers for wave equation with momentum damping in first-order formulation.}
	\label{fig:wave}
\end{figure}
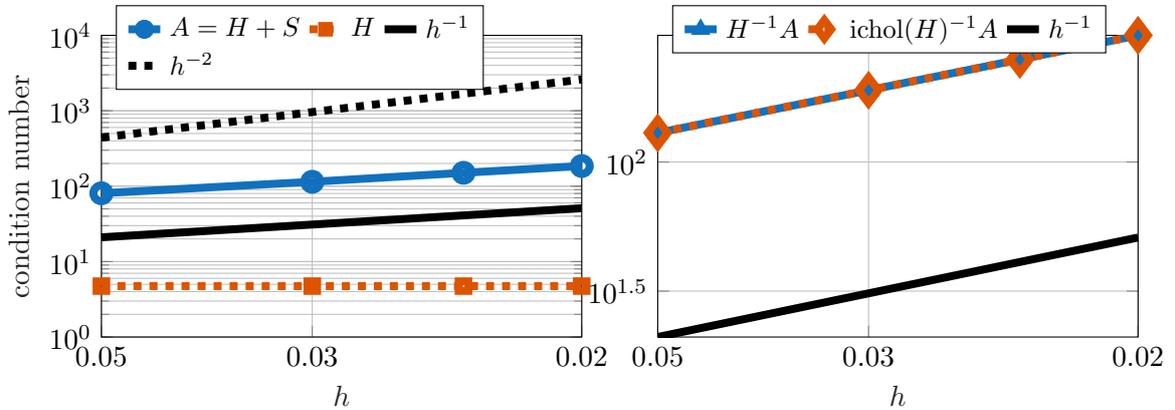

 We illustrate the unboundedness  in Proposition~\ref{prop:wave} in Figure~\ref{fig:wave}, where we choose $\Omega = [0,1]$ and $\rho=\nu =1$. In the left plot, we see that the condition number of $A$ and also that of the skew-symmetric part consisting of the first-order differential operators in the off-diagonal blocks behaves like $\mathcal{O}(h^{-1})$. In the right plot we observe that the preconditioning has no effect, which is to be expected as $\calH = I \in L(L^2(\Omega),L^2(\Omega))$.

In this subsection we have demonstrated that the preconditioning with the symmetric part does not lead to a mesh-independent order of the condition number if the skew-symmetric part is dominant in the sense that it has a higher order of the differential operator in the skew-symmetric part.

In this example, it would hence clearly be more suitable to choose the (unbounded) skew-symmetric part $\calS$ or an approximation thereof as a preconditioner. In this regard, we refer to the work \cite{golub2000preconditioning} focusing on preconditioning with the skew-symmetric part. An alternative would be to multiply the equation with the imaginary unit such that $i\calS$ self-adjoint. 

\subsection{Beam equation with structural (Kelvin-Voigt) damping}\label{subsec:kv}

As a last example, we consider a beam equation that is subject to strong (structural) damping on the domain $\Omega = (0,1)$.  %
Denoting by $x:\R_{\geq 0}\times (0,1) \to \R$ the transverse displacement of the beam, the dynamics of the beam is given by
the PDE
\begin{equation*}
    \ddot x(t,r) + \frac{\partial^2}{\partial r^2}\left(E \frac{\partial ^2}{\partial r^2}x(t,r) + C\frac{\partial^3}{\partial r^2\partial t}x(t,r)\right) = f(t,r),
\end{equation*}
with boundary conditions
\begin{equation}\label{eq:bcs_beam}
        x(t,0)=\frac{\partial}{\partial r}x(t,1) = \frac{\partial^2}{\partial r^2}x(t,0) = \frac{\partial^3}{\partial r^3}x(t,1)=0.
\end{equation}
For an in-depth analytical treatment, we refer to \cite{jacob2008analyticity}. Note that although at first sight, the beam equation seems hyperbolic, it was proven in \cite{jacob2008analyticity} that it gives rise to an analytic semigroup, hence behaving more like a parabolic equation such as the advection-diffusion-reaction equation. The reason for this is the strong damping. This parabolic nature will also be observed in the boundedness of the preconditioned operator and hence uniform condition numbers of Galerkin projections.

To obtain a first order dissipative %
formulation, we define the variables $p = \dot x$, $q = \frac{\partial^2}{\partial r^2}x$ which yield the dynamical system 
\begin{equation*}
\frac{\mathrm{d}}{\mathrm{d}t}
    \begin{bmatrix}
        p\\
        q
    \end{bmatrix}
    = \begin{bmatrix}
        -\frac{\partial^4}{\partial r^4} & -\frac{\partial^2}{\partial r^2}\\
       \frac{\partial^2}{\partial r^2} & 0
    \end{bmatrix}
    \begin{bmatrix}
        p\\q
    \end{bmatrix}
    + \begin{bmatrix}
        f\\
        0
    \end{bmatrix}.
\end{equation*}
As the boundary conditions for the position \eqref{eq:bcs_beam} also translate to the momentum, we have for $p_1,p_2:(0,1)\to \R$ satisfying the boundary conditions in \eqref{eq:bcs_beam} the (formal) equality
\begin{align*}
    &\langle p_1^{(4)}, p_2\rangle_{L^2(0,1)} \\&= \langle p_1^{(2)},p_2^{(2)}\rangle_{L^2(0,1)} + p_1^{(3)}(1)p_2(1) - p_1^{(3)}(0)p_2(0) - p_1^{(2)}(1)p_2^{(1)}(1) + p_1^{(2)}(0)p_2^{(1)}(0).
     \\&= \langle p_1^{(2)},p_2^{(2)}\rangle_{L^2(0,1)} 
\end{align*}
and, for momentum $p$ and curvature $q$ respecting the boundary conditions 
\eqref{eq:bcs_beam}, we have
\begin{align*}
    \langle p^{(2)},q \rangle_{L^2(0,1)} &= \langle p,q^{(2)}\rangle_{L^2(0,1)} + p^{(1)}(1) q(1) - p^{(1)}(0)q(0) - p(1)q^{(1)}(1) + p(0)q^{(1)}(0) \\&= \langle p,q^{(2)}\rangle_{L^2(0,1)}. 
\end{align*}
Hence, the corresponding operator in weak form (readily accessible for finite element methods) may be formulated via
\begin{align}\label{KVdampingop}
	\calA = \begin{bmatrix}
		\calA_{11} & \calA_{12} \\
		\calA_{21} & 0
	\end{bmatrix}
\end{align}
where, setting $\dom(\calA_{11}) := \{p\in H^2(0,1)\,|\, p(0)=p'(1)=0\}\subset H^2(0,1)$ and $H^{-2}(0,1) := \dom(\calA_{11})^*$, where the dual space is taken w.r.t.\ the Pivot space $L^2(0,1)$, we consider the operators
\begin{align*}
    \calA_{11}: H^{-2}(0,1) \supset \dom(\calA_{11}) \to H^{-2}(0,1)&, \quad 
    \calA_{12}: L^2(0,1) \to H^{-2}(0,1),\\
    \calA_{21}&: \dom(\calA_{11}) \to L^2(0,1)
\end{align*}
which are for $p_1,p_2\in \dom(\calA_{11})$ and $q\in L^2(0,1)$ defined by
\begin{align*}
    \langle\calA_{11}p_1,p_2\rangle_{H^{-2}(0,1), H^2(0,1)} = \langle p_1^{(2)},p_2^{(2)}\rangle_{L^2(\Omega)}&, \ \langle\calA_{12}q,p_1\rangle_{H^{-2}(0,1), H^2(0,1)}:= \langle q,p_1^{(2)}\rangle_{L^2(0,1)}\\
    \langle \calA_{21} p_1,q\rangle_{L^2(0,1)} &= -\langle p_1^{(2)} ,q\rangle_{L^2(0,1)}.
\end{align*}
As in the Oseen example discussed in Subsection~\ref{subsec:oseen},  we form the Schur complement
\begin{equation}\label{eq:schurkv}
    \calW := -\calA_{21} \calA_{11}^{-1} \calA_{12},
\end{equation}
and we obtain an analogous result to Proposition~\ref{prop:oseen}.
\begin{proposition}\label{prop:kelvin}
Consider the operator $\calA$ in \eqref{KVdampingop} and the associated Schur complement $\calW$ in~\eqref{eq:schurkv}.  Then $\calW$ is bounded, i.e., $\calW\in L(L^2(0,1), L^2(0,1))$, and accretive.
\end{proposition}
\begin{proof}
    The proof is analogous to the proof of Proposition~\ref{prop:oseen}. First, we show that $\calA_{11}:H^{-2}(0,1) \supset \dom(\calA_{11}) \to H^{-2}(0,1)$ gives rise to a bounded and bilinear form $a:\dom(\calA_{11}) \times \dom(\calA_{11}) \to \R$ defined by
    \begin{equation*}
        a(p_1,p_2) := \langle \calA_{11} p_1,p_2\rangle_{H^{-2}(0,1),H^2(0,1)}
    \end{equation*}
    to invoke the Lax Milgram theorem. To this end, for $p_1,p_2 \in \dom(\calA_{11})$, boundedness follows by
    \begin{equation*}
        a(p_1,p_2) = \langle p_1^{(2)},p_2^{(2)}\rangle_{L^2(0,1)} \leq \|p_1\|_{H^2(0,1)} \|p_2\|_{H^2(0,1)}
    \end{equation*}
    from the Cauchy-Schwarz inequality. Moreover, for $p\in \dom(\calA_{11})$,
    \begin{align*}
        a(p,p) = \|p^{(2)}\|^2_{L^2(0,1)} &= \frac13 \|p^{(2)}\|^2_{L^2(0,1)} + \frac23 \|p^{(2)}\|^2_{L^2(0,1)} \\
        &\geq \frac13 \|p^{(2)}\|^2_{L^2(0,1)} + \frac23 c\|p^{(1)}\|^2_{L^2(0,1)}
    \end{align*}
    by applying the Poincaré inequality to $p^{(1)}$ in view of the boundary condition $p^{(1)}(1)=0$ for $p \in \dom(\calA_{11})$. Again applying the Poincaré inequality for $p$, using that $p(0)=0$ for $p \in \dom(\calA_{11})$, we deduce
    \begin{equation*}
        a(p,p) \geq \frac13 \left(\|p^{(2)}\|^2_{L^2(0,1)} + c \|p^{(1)}\|^2_{L^2(0,1)} + c^2\|p\|^2_{L^2(0,1)} \right) = \tilde c\|p\|^2_{H^2(0,1)}
    \end{equation*}
    with $\tilde c = \frac13\min \{1,c^2\}$. Thus, by the Lax Milgram theorem $\calA_{11}: \dom(\calA_{11}) \to H^{-2}(0,1)$ has a bounded inverse $\calA_{11}^{-1}\in L(H^{-2}(0,1),H^2(0,1))$ with $\operatorname{im}\calA_{11}^{-1}\subset \operatorname{dom}(\calA_{11})$. Consequently, as a concatenation of bounded linear maps, the Schur complement $\calW$ defined in \eqref{eq:schurkv} is bounded, as well. 

    By the above computations it follows also that $\calA_{11}$ is coercive. Moreover, $\calA_{21}=-\calA_{12}^*$ follows by definition and integration by parts using the boundary conditions included in $\dom(\calA_{11})$. Consequently, by an analogous argumentation as in the proof of Proposition~\ref{prop:oseen}, the Schur complement is dissipative.
\end{proof}

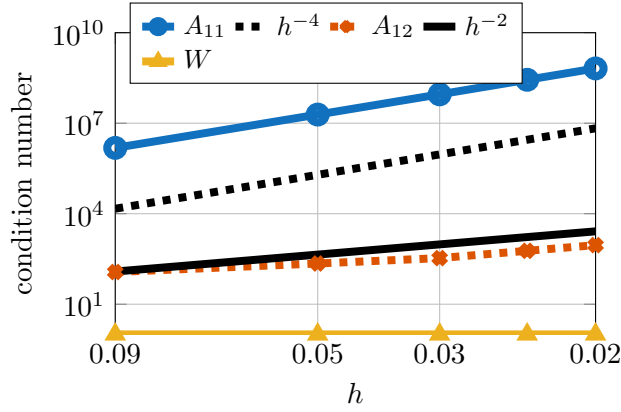
\begin{figure}[htb]
	\centering
\definecolor{mycolor1}{rgb}{0.06600,0.44300,0.74500}%
\definecolor{mycolor2}{rgb}{0.86600,0.32900,0.00000}%
\definecolor{mycolor3}{rgb}{0.92900,0.69400,0.12500}%
\definecolor{mycolor4}{rgb}{0.12941,0.12941,0.12941}%

\begin{tikzpicture}

\begin{axis}[%
  width=.4\linewidth,
  height=4cm,
  at={(0in,0.555in)},
  scale only axis,
  x dir=reverse,
  xmode=log,
  xmin=0.0196078431372549,
  xmax=0.0909090909090909,
  xminorticks=true,
  xlabel={$h$},
  ymode=log,
  ymin=1,
  ymax=10000000000,
  yminorticks=true,
  ylabel={condition number},
  axis background/.style={fill=white},
  xmajorgrids, xminorgrids,
  ymajorgrids, yminorgrids,
  legend style={at={(0.03,1.1)}, anchor=north west, legend columns=4, legend cell align=left, align=left},
  legend image post style={scale=0.7},
  legend style={font=\small},
  xtick={0.166666666666667,0.0909090909090909,0.0476190476190476,0.032258064516129,0.0196078431372549},
  xticklabels={0.17,0.09,0.05,0.03,0.02},
]

\addplot [color=mycolor1, line width=3.0pt, mark size=3.0pt, mark=o, mark options={solid, mycolor1}]
  table[row sep=crcr]{%
0.0909090909090909  1518019.55556606\\
0.0476190476190476  19330926.2185746\\
0.032258064516129   90416721.8808671\\
0.024390243902439   274508738.513109\\
0.0196078431372549  654095859.389819\\
};
\addlegendentry{$A_{11}$}

\addplot [color=black, dashed, line width=3.0pt]
  table[row sep=crcr]{%
0.0909090909090909  14641\\
0.0476190476190476  194481\\
0.032258064516129   923521\\
0.024390243902439   2825761\\
0.0196078431372549  6765201\\
};
\addlegendentry{$h^{-4}$}

\addplot [color=mycolor2, dashed, line width=3.0pt, mark size=3.0pt, mark=x, mark options={solid, mycolor2}]
  table[row sep=crcr]{%
0.0909090909090909  116.628370938722\\
0.0476190476190476  223.620052224127\\
0.032258064516129   337.58882774891\\
0.024390243902439   588.323129370522\\
0.0196078431372549  908.226791330851\\
};
\addlegendentry{$A_{12}$}

\addplot [color=black, line width=3.0pt]
  table[row sep=crcr]{%
0.0909090909090909  121\\
0.0476190476190476  441\\
0.032258064516129   961\\
0.024390243902439   1681\\
0.0196078431372549  2601\\
};
\addlegendentry{$h^{-2}$}

\addplot [color=mycolor3, line width=3.0pt, mark size=1.7pt, mark=triangle, mark options={solid, mycolor3}]
  table[row sep=crcr]{%
0.0909090909090909  1.00000000001058\\
0.0476190476190476  1.00000000016005\\
0.032258064516129   1.00000000204369\\
0.024390243902439   1.00000000185927\\
0.0196078431372549  1.00000000478196\\
};
\addlegendentry{$W$}

\end{axis}

\end{tikzpicture}%
 	\caption{Condition numbers for beam equation with Kelvin-Voigt damping.}
	\label{fig:kv}
\end{figure}

 The numerical results for this example are presented in Figure~\ref{fig:kv}. While the condition number of the top left block of the operator behaves like $\mathcal{O}(h^{-4})$ due to the presence of four derivatives and the condition number of the top right block scales like $\mathcal{O}(h^{-2})$ due to the presence of two derivatives, the condition number of the Schur complement is bounded uniformly in the mesh width. This reflects the boundedness proven in Proposition~\ref{prop:kelvin}. 

So far we have studied operator preconditioning for linear operators which have a symmetric part that is semidefinite. In the next section we apply these techniques in the solution of optimal control problems.
\section{Application in optimal control}\label{sec:optcont}
Having discussed the condition numbers in the context of operator preconditioning for partial differential equations, in this section we  discuss two implementations of the structure-exploiting iterative schemes in the context of optimal control problems. We consider the prototypical problem
\begin{equation}\label{eq:OCP}
    \min_{(x,u)\in \dom(\mathcal{A})\times U} \frac12\|\mathcal{C}x-y_\mathrm{ref}\|_Y^2 + \frac\lambda2 \|u-u_\mathrm{ref}\|_U^2 \quad \mathrm{s.t.} \quad \mathcal{A}x-\mathcal{B}u = f,\ y=\mathcal C x
\end{equation}
for Hilbert spaces $X,Y,U$. Here, $\mathcal{A}:\dom(\mathcal{A})\subset X\to X$ is a densely defined linear operator, $\mathcal{B}\in L(U,X)$ is an input operator, $\mathcal{C}\in L(X,Y)$ is an output operator, $\lambda >0$ is a regularization parameter, $f\in X$ is a source term and $y_\mathrm{ref}\in Y$ as well as $u_\mathrm{ref}\in U$ are a reference output and reference control.

Throughout the following we assume that $\mathcal{A}$ is an accretive operator giving rise to a splitting 
\begin{equation}\label{eq:structureagain}
    \calA = \calH + \calS
\end{equation}
with symmetric $\calH$ and skew-symmetric $\calS$. For particular applications in stationary problems, we refer to the examples provided in Subsections~\ref{subsec:advec}--\ref{subsec:kv}. 
Note, that $\mathcal{A}$ may be the model of a semi-discretized time-dependent problem, as long as the time discretization is dissipativity-preserving, as discussed in Section after~\eqref{eq:timedep}. Using, e.g., an implicit midpoint discretization of a dissipative evolution equation $\dot x = -\calM x$ with dissipative operator $\calM$ as in \eqref{eq:timedep},
with time step size $\delta t > 0$, the operator $\mathcal{A}$ in the optimal control problem \eqref{eq:OCP} reads
\begin{align*}
    \calA = 
\left(\begin{smallmatrix}
I & 0 & 0 & 0 & \cdots & 0 \\
I-\nicefrac{\delta t}{2}\calM & I+\nicefrac{\delta t}{2}\calM  & 0 & 0 & \cdots & 0 \\
0 & I-\nicefrac{\delta t}{2}\calM &I+\nicefrac{\delta t}{2}\calM  & 0 & \cdots & 0 \\
\vdots & \vdots & \ddots & \ddots & \ddots & \vdots \\
\vdots & \vdots & \vdots & I-\nicefrac{\delta t}{2}\calM & I+\nicefrac{\delta t}{2}\calM  & 0 \\
0 & 0 & 0 & 0 & I-\nicefrac{\delta t}{2}\calM & I+\nicefrac{\delta t}{2}\calM  \\
\end{smallmatrix}\right)
\end{align*}
where the first block line assigns a given initial condition. Denoting by $\calH_\calM$ and $\calS_\calM$ the symmetric resp.\ skew-symmetric part of $\calM$, we may split
\[
I+\nicefrac{\delta t}{2}\mathcal{M} = \underbrace{I + \nicefrac{\delta t}{2}\calH_{\calM}}_{:=\calH_d} \underbrace{- \nicefrac{\delta t}{2}\calS_{\calM}}_{:= \calS_d}.
\]
Clearly, $\calH_d$ is symmetric positive definite and $\calS_d$ is skew-symmetric. Consequently, we set
\begin{equation*}
    \calA = \calH + \calS, \qquad \calH = \calH^* \geq I, \quad \calS = -\calS^*
\end{equation*}
with
\begin{equation*}
   \calH = \frac12\left(\begin{smallmatrix}
I & I-\nicefrac{\delta t}{2}\calM^* & 0 & 0 & \cdots &  \\
I-\nicefrac{\delta t}{2}\calM & 2\calH_d  & I-\nicefrac{\delta t}{2}\calM^* & 0 & \cdots & 0 \\
0 & I-\nicefrac{\delta t}{2}\calM &2\calH_d  & I-\nicefrac{\delta t}{2}\calM^* & \cdots & 0 \\
\vdots & \vdots & \ddots & \ddots & \ddots & \vdots \\
\vdots & \vdots & \vdots & I-\nicefrac{\delta t}{2}\calM & 2\calH_d  & I-\nicefrac{\delta t}{2}\calM^*\\
0 & 0 & 0 & 0 & I-\nicefrac{\delta t}{2}\calM & 2\calH_d \\
\end{smallmatrix}\right)
\end{equation*}
and
\begin{equation*}
   \calS = \frac12\left(\begin{smallmatrix}
0 & -I+\nicefrac{\delta t}{2}\calM^* & 0 & 0 & \cdots &  \\
I-\nicefrac{\delta t}{2}\calM & 2\calS_d  & -I+\nicefrac{\delta t}{2}\calM^* & 0 & \cdots & 0 \\
0 & I-\nicefrac{\delta t}{2}\calM &2\calS_d & -I+\nicefrac{\delta t}{2}\calM^* & \cdots & 0 \\
\vdots & \vdots & \ddots & \ddots & \ddots & \vdots \\
\vdots & \vdots & \vdots & I-\nicefrac{\delta t}{2}\calM & 2\calS_d  & -I+\nicefrac{\delta t}{2}\calM^*\\
0 & 0 & 0 & 0 & I-\nicefrac{\delta t}{2}\calM & 2\calS_d \\
\end{smallmatrix}\right).
\end{equation*}
Note that one would usually not apply a factorization-based solver to the full block operator $\calA$, but rather solve the resulting system block-wise recovering the iteration \eqref{eq:IMP}.

\begin{remark}
    In principle, we could also consider time-dependent problems in function space. There, $X = L^2(0,T;\bar X)$ for another Hilbert space $\bar X$, $\calA = \tfrac{\mathrm{d}}{\mathrm{d}t} - \calM$ with dissipative $\calM:\dom(\calM)\subset \bar X \to \bar X$ being the generator of a strongly continuous semigroup, hence closed such that $\dom(\calM)$ becomes a Banach space when endowed with the graph norm. Correspondingly, we may set $\dom(\calA) = H^1(0,T;\bar X)\cap L^2(0,T;\dom(\calM))$ such that $\calA$ is densely defined in $L^2(0,T;\bar X)$, see e.g.\ \cite{farkas2025dissipativity}. For parabolic equations where $\calM$ gives rise to an analytic semigroup as the ones considered Section~\ref{sec:conditioning} with exception of the wave equation, $\calA$ is closed due to the maximal parabolic regularity, see~\cite[Sec.\ 3.6]{BensDaPr07}. For example, if $-\calM$ is a second order elliptic operator on $\bar X=H^{-1}(\Omega)$ in weak form, $D(\calA) = H^1(0,T;H^{-1}(\Omega))\cap L^2(0,T;H^1(\Omega))$ corresponds to the usual $W(0,T)$-space used in variational theory~\cite{Tro10}.
\end{remark}

\medskip

\textbf{Optimality conditions.} Before going to the structure-exploiting linear solvers, we briefly recall the optimality conditions for \eqref{eq:OCP}, assuming that $\calA$ has closed range, see~\cite[Theorem 1.1 and Remark 1.2]{Schiela13}. Note that this closed range assumption is satisfied e.g.\ if $\calA$ is surjective which is the case for differential operators in a suitable functional analytic setting (see e.g.\ the operators in Subsection~\ref{subsec:advec}--\ref{subsec:kv}). 

Let $(x,u)\in \dom(\calA)\times {U}$ be optimal. Then, there is $p\in \dom(\calA^*)$ such that 
\begin{equation}\label{eq:optsys}
\begin{bmatrix}
    \calC^*\calC & 0 & \calA^*\\
    0 & \lambda I & -\calB^*\\
    \calA& -\calB &0
\end{bmatrix}\begin{bmatrix}
    x\\u\\p
\end{bmatrix}
= 
\begin{bmatrix}
    \calC^*y_\mathrm{ref}\\
    \lambda\, u_\mathrm{ref}\\
    f
\end{bmatrix}.
\end{equation}

We now discuss three different routes to leverage symmetric/skew-symmetric splittings in the optimality system \eqref{eq:optsys}. 

\subsection{Direct splitting of the optimality system}
In a first approach, we do not require any structure of the constraint operator $\calA$. Multiplying the last block row in \eqref{eq:optsys} by minus one, we may split the optimality system into its symmetric and skew-symmetric part via
\begin{equation*}
    \calH =     \begin{bmatrix}
        \calC^*\calC & 0 & 0\\
        0 & \lambda I & 0\\
        0& 0 &0
        \end{bmatrix}, \qquad \calS = \begin{bmatrix}
        0 & 0 & \calA^*\\
        0 & 0 & -\calB^*\\
        -\calA& \calB &0
        \end{bmatrix}.
\end{equation*}
Here, the symmetric part is not invertible. If  $\calC^*\calC$ is invertible, then we may apply the Schur complement to the top left 2x2 block of the symmetric part. This leads to a problem governed by the matrix %
\begin{equation*}
    \begin{bmatrix}
        \calA & \calB
    \end{bmatrix} \begin{bmatrix}
        \calC^*\calC & 0\\
        0&\lambda I
    \end{bmatrix}^{-1}
    \begin{bmatrix}
        \calA^*\\
        \calB^*
    \end{bmatrix} = \calA (\calC^*\calC)^{-1} \calA^* + \frac1\lambda \calB\calB^*.
\end{equation*}
This operator is (formally) symmetric and semidefinite. If $\calA$ is an isomorphism, as in many PDE applications, the first term is even definite. Hence, the method of choice would be a suitable variant of the conjugate gradient method.  Note, however, that the invertibility of $\calC^*\calC$ which, in many applications with partial differential equations, implies observation on the whole domain would be very restrictive. This problem could be resolved by a further partitioning of $ \mathcal C$ in the observable and unobservable part and then using another Schur complement, again an approach that is often not feasible for large scale problems. This formulation is often used for designing preconditioners, typically for the case when $\calC^*\calC$ is invertible, see \cite{pearson2012regularization,rees2010optimal,schoberl2007symmetric} among others. In case of non-invertible observation term, for the sake of preconditioning a perturbation rendering $\calC^*\calC$ invertible could be used but still would only be useful when dealing with the full saddle point formulation.

\subsection{Condensed formulation}\label{sec:reduced}
    A different but also well-known approach in optimal control of partial differential equations is the reduction to the control variable via the state-to-control map, see e.g.\ \cite{Tro10}. For this, we assume that the operator $\calA$ in the optimal control problem \eqref{eq:OCP} is boundedly invertible such that we may eliminate the state  via $x = \calA^{-1}\calB u + \calA^{-1}f$. This leads to the reduced unconstrained optimization problem
    \begin{equation*}
        \min_{u\in U} f(u) := \frac12\|\calC(\calA^{-1}\calB u + \calA^{-1}f) - y_\mathrm{ref}\|_Y^2 + \frac\lambda2 \|u-u_\mathrm{ref}\|_U^2.
    \end{equation*}
    Due to the strict convexity of the cost function,
    solving this problem amounts to solving the first-order necessary (and sufficient) optimality condition at the optimal control $u^*\in U$
    \begin{equation}\label{eq:cg}
        0 = \nabla f(u^*) = \left((\calC\calA^{-1}\calB)^*\calC\calA^{-1}\calB+ \lambda I\right)u^* + (\calC\calA^{-1}\calB)^*\calC(\calA^{-1}f - y_\mathrm{ref}).
    \end{equation}
    Here, the governing operator
    \begin{equation}\label{eq:symmetricmatrix}
    (\calC\calA^{-1}\calB)^*\calC\calA^{-1}\calB+ \lambda I
    \end{equation}
    is symmetric and strictly accretive such that we can apply the conjugate gradient method. However, an application of this matrix requires in each step the solution of the state and the adjoint equation, i.e., evaluating 
    $\calA^{-1}$ and $\calA^{-*}$.

    To implement these solution operators, we use the  accretivity and the symmetric/skew-symmetric splitting of the governing operator $\calA = \calH + \calS$, which clearly also carries over to the adjoint operator $\calA^* = \calH^* + \calS^*$. More precisely, we will solve both the state and the adjoint equation via GMRES, Rapoport's or Widlund's method endowed with preconditioning by the symmetric part or an approximation of it
    such as an incomplete Cholesky factorizations or an algebraic multigrid method. While we ignore a possible source of inexactness in these approximations, we would like to again mention the work \cite{diab2023flexible}, where flexible methods that tolerate inexact evaluations of the preconditioner are discussed.

    In the following we denote again by $A$, $H$ and $S$ the Galerkin projections of the operators $\calA$, $\calH$ and $\calS$ onto a suitable finite element space, respectively. We then use the following methods and their abbreviations. 
    
     \textbf{Methods for solving $\calA^{-1}$ and $\calA^{-*}$.} The methods considered in the following to evaluate state and adjoint equation in the outer CG loop applied to the system governed by the (discretization of the) operator \eqref{eq:symmetricmatrix} are:
    \begin{itemize}[leftmargin=2em]
        \item Direct: Direct solution by means of an LU factorization of $A$.
        \item ILU: Incomplete LU factorization of $A$ with prescribed drop tolerance.
        \item GMRES: Generalized Minimal Residual Method for $Ax=b$ and $A^*p = b$ without preconditioner.
        \item GMRES(IC): GMRES with an incomplete Cholesky factorization $LL^\top \approx H$ as left preconditioner.
        \item GMRES(MG): GMRES with algebraic multigrid method for $H$ as a left preconditioner. Here, we use the algebraic multigrid method as part of the HSL library~\cite{hsl2023}.
        \item Rapoport: Rapoport's method~\cite{rapoport1978nonlinear} as described in Section~\ref{sec:solvers}.
        \item Widulund: Widlund's method~\cite{widlund1978lanczos} as described in Section~\ref{sec:solvers}.
    \end{itemize}
    
We evaluate the suggested approach by means of two applications, namely the advection-diffusion-reaction equation from Subsection~\ref{subsec:advec} and the Stokes equation from Subsection~\ref{subsec:stokes}. For both examples, we solve the reduced optimality system \eqref{eq:cg} with a conjugate gradient method up to relative tolerance $\mathrm{cgtol}>0$. In each step, we solve state and adjoint equation with settings described in Table~\ref{tab:methods_params}. Note that for inner solvers requiring a stopping criterion, we use the inner relative tolerance $\mathrm{cgtol}/10$ adapted to the tolerance of the outer cg loop.

\begin{table}[htb]
\centering
\scalebox{.75}{
\begin{tabular}{c | c c c c c c c}
\hline
& direct & ILU & GMRES & GMRES(IC)  & GMRES(MG)  & Rapoport & Widlund  \\
\hline
assembly & lu(A) & ilu(A,$\mathrm{cgtol}$/100) & - &  ichol(H,$10^{-1}$)& multigrid (2c) & multigrid (2c) & multigrid (2c) \\
tolerance & - &- & $\mathrm{cgtol}/10$ & $\mathrm{cgtol}/10$ & $\mathrm{cgtol}/10$  & $\mathrm{cgtol}/10$ &$\mathrm{cgtol}/10$  \\
\hline
\end{tabular}
}
\caption{Used inner solver for the outer CG iteration with tolerance $\mathrm{cgtol}$. \label{tab:methods_params}}
\end{table}

\subsubsection{Advection-diffusion-reaction equation}
As first example,  consider a three-dim\-ensional stationary advection-diffusion-reaction equation introduced in Subsection~\ref{subsec:advec} on the unit cube $\Omega = [0,1]^3$. We choose constant diffusivity $\nu \equiv 1$, advection term $\textbf{b}\equiv \begin{bmatrix} -0.5&0&0\end{bmatrix}^\top$, reaction term $\textbf{c}\equiv 1$ and a source term $f\equiv 10$. For the input-output configuration, we choose full observation and control, that is, $X=Y=U=L^2(\Omega)$, $B= C = I$ and vanishing reference functions $y_\mathrm{ref} = u_\mathrm{ref} \equiv 0$.

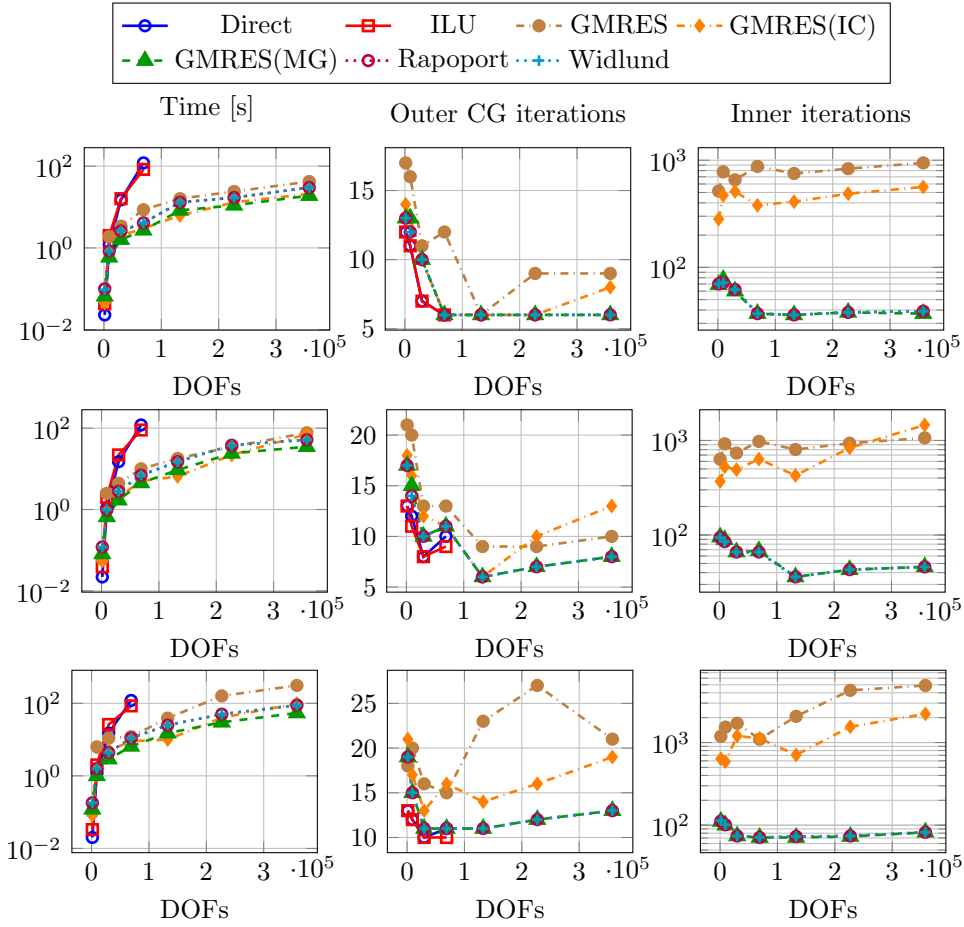
\begin{figure}[htb]
\centering
\scalebox{1}{%

\begin{subfigure}[t]{136.9pt}
\centering
\begin{tikzpicture}
\begin{axis}[
    title={Time [s]},
    xlabel={DOFs},
    ymode=log,
    width=\linewidth,
    height=4cm,
    grid=both,
    tick label style={font=\small},
    label style={font=\small},
    title style={font=\small},
    cycle list name=mycyclelist,
    legend style={font=\small,at={(1.7,1.35)}, anchor=south, align=center},
    legend columns = 4,
    every x tick scale label/.style={xshift = 2.9cm,yshift=-1.3cm}, %
]
\addplot coordinates {(1331,0.023) (9261,1.2) (29791,15) (68921,120)};
\addlegendentry{Direct}
\addplot coordinates {(1331,0.043) (9261,2) (29791,16) (68921,83)};
\addlegendentry{ILU}
\addplot coordinates {(1331,0.052) (9261,1.9) (29791,3.4) (68921,8.6) (132651,16) (226981,24) (357911,42)};
\addlegendentry{GMRES}
\addplot coordinates {(1331,0.045) (9261,0.86) (29791,2) (68921,3) (132651,6.2) (226981,13) (357911,21)};
\addlegendentry{GMRES(IC)}
\addplot coordinates {(1331,0.066) (9261,0.59) (29791,1.6) (68921,2.7) (132651,8.3) (226981,11) (357911,19)};
\addlegendentry{GMRES(MG)}
\addplot coordinates {(1331,0.1) (9261,0.84) (29791,2.6) (68921,4.1) (132651,13) (226981,17) (357911,30)};
\addlegendentry{Rapoport}
\addplot coordinates {(1331,0.095) (9261,0.84) (29791,2.6) (68921,4.1) (132651,13) (226981,17) (357911,30)};
\addlegendentry{Widlund}
\end{axis}
\end{tikzpicture}
\end{subfigure}\hspace{-.9cm}
\begin{subfigure}[t]{136.9pt}
\centering
\begin{tikzpicture}
\begin{axis}[
    title={Outer CG iterations },
    xlabel={DOFs},
    width=\linewidth,
    height=4cm,
    grid=both,
    tick label style={font=\small},
    label style={font=\small},
    title style={font=\small},
    cycle list name=mycyclelist,
    every x tick scale label/.style={xshift = 2.9cm,yshift=-.4cm}, %
]
\addplot coordinates {(1331,12) (9261,11) (29791,7) (68921,6)};
\addplot coordinates {(1331,12) (9261,11) (29791,7) (68921,6)};
\addplot coordinates {(1331,17) (9261,16) (29791,11) (68921,12) (132651,6) (226981,9) (357911,9)};
\addplot coordinates {(1331,14) (9261,13) (29791,10) (68921,6) (132651,6) (226981,6) (357911,8)};
\addplot coordinates {(1331,13) (9261,13) (29791,10) (68921,6) (132651,6) (226981,6) (357911,6)};
\addplot coordinates {(1331,13) (9261,12) (29791,10) (68921,6) (132651,6) (226981,6) (357911,6)};
\addplot coordinates {(1331,13) (9261,12) (29791,10) (68921,6) (132651,6) (226981,6) (357911,6)};
\end{axis}
\end{tikzpicture}
\end{subfigure}\hspace{-.9cm}
\begin{subfigure}[t]{136.9pt}
\centering
\begin{tikzpicture}
\begin{axis}[
    title={Inner iterations},
    xlabel={DOFs},
    ymode=log,
    width=\linewidth,
    height=4cm,
    grid=both,
    tick label style={font=\small},
    label style={font=\small},
    title style={font=\small},
    legend style={font=\tiny},
    legend pos=north west,
    cycle list name=mycyclelist,
    every x tick scale label/.style={xshift = 2.9cm,yshift=1.8cm}, %
]
\addplot coordinates {};
\addplot coordinates {};
\addplot coordinates {(1331,515) (9261,779) (29791,656) (68921,875) (132651,752) (226981,834) (357911,944)};
\addplot coordinates {(1331,282) (9261,469) (29791,511) (68921,379) (132651,410) (226981,487) (357911,563)};
\addplot coordinates {(1331,69) (9261,77) (29791,60) (68921,37) (132651,36) (226981,38) (357911,37)};
\addplot coordinates {(1331,70) (9261,72) (29791,62) (68921,37) (132651,36) (226981,38) (357911,39)};
\addplot coordinates {(1331,70) (9261,72) (29791,62) (68921,37) (132651,36) (226981,38) (357911,39)};
\end{axis}
\end{tikzpicture}
\end{subfigure}
\hfill  \vspace*{.2cm}

}\\%
\scalebox{1}{%
\hspace{-.25cm}
\begin{subfigure}[t]{136.9pt}
\centering
\begin{tikzpicture}
\begin{axis}[
    xlabel={DOFs},
    ymode=log,
    width=\linewidth,
    height=4cm,
    grid=both,
    tick label style={font=\small},
    label style={font=\small},
    title style={font=\small},
    cycle list name=mycyclelist,
    every x tick scale label/.style={xshift = 2.9cm,yshift=-1.3cm}, %
]
\addplot coordinates {(1331,0.022) (9261,1.1) (29791,15) (68921,120)};
\addplot coordinates {(1331,0.039) (9261,2) (29791,22) (68921,90)};
\addplot coordinates {(1331,0.061) (9261,2.5) (29791,4.4) (68921,10) (132651,18) (226981,35) (357911,77)};
\addplot coordinates {(1331,0.054) (9261,0.98) (29791,2) (68921,5) (132651,6.5) (226981,22) (357911,71)};
\addplot coordinates {(1331,0.081) (9261,0.66) (29791,1.7) (68921,4.5) (132651,9.4) (226981,24) (357911,35)};
\addplot coordinates {(1331,0.12) (9261,0.98) (29791,2.8) (68921,7) (132651,15) (226981,38) (357911,52)};
\addplot coordinates {(1331,0.12) (9261,1) (29791,2.8) (68921,7) (132651,15) (226981,38) (357911,52)};
\end{axis}
\end{tikzpicture}
\end{subfigure}\hspace{-.74cm}
\begin{subfigure}[t]{136.9pt}
\centering
\begin{tikzpicture}
\begin{axis}[
    xlabel={DOFs},
    width=\linewidth,
    height=4cm,
    grid=both,
    tick label style={font=\small},
    label style={font=\small},
    title style={font=\small},
    cycle list name=mycyclelist,
    legend style={font=\tiny, at={(1,1.2)}, anchor=north east},
    legend columns=1,
    every x tick scale label/.style={xshift = 2.9cm,yshift=-.4cm}, %
]
\addplot coordinates {(1331,13) (9261,12) (29791,8) (68921,10)};
\addplot coordinates {(1331,13) (9261,11) (29791,8) (68921,9)};
\addplot coordinates {(1331,21) (9261,20) (29791,13) (68921,13) (132651,9) (226981,9) (357911,10)};
\addplot coordinates {(1331,18) (9261,16) (29791,12) (68921,11) (132651,6) (226981,10) (357911,13)};
\addplot coordinates {(1331,17) (9261,15) (29791,10) (68921,11) (132651,6) (226981,7) (357911,8)};
\addplot coordinates {(1331,17) (9261,14) (29791,10) (68921,11) (132651,6) (226981,7) (357911,8)};
\addplot coordinates {(1331,17) (9261,14) (29791,10) (68921,11) (132651,6) (226981,7) (357911,8)};
\end{axis}
\end{tikzpicture}
\end{subfigure}\hspace{-.9cm}
\begin{subfigure}[t]{136.9pt}
\centering
\begin{tikzpicture}
\begin{axis}[
    xlabel={DOFs},
    ymode=log,
    width=\linewidth,
    height=4cm,
    grid=both,
    tick label style={font=\small},
    label style={font=\small},
    title style={font=\small},
    cycle list name=mycyclelist,
    legend style={font=\tiny},
    legend pos=north west,
    every x tick scale label/.style={xshift = 2.9cm,yshift=1.53cm}, %
]
\addplot coordinates {};
\addplot coordinates {};
\addplot coordinates {(1331,633) (9261,921) (29791,735) (68921,977) (132651,801) (226981,934) (357911,1063)};
\addplot coordinates {(1331,369) (9261,526) (29791,491) (68921,634) (132651,426) (226981,834) (357911,1458)};
\addplot coordinates {(1331,95) (9261,89) (29791,66) (68921,68) (132651,36) (226981,43) (357911,46)};
\addplot coordinates {(1331,94) (9261,85) (29791,66) (68921,66) (132651,36) (226981,43) (357911,46)};
\addplot coordinates {(1331,94) (9261,85) (29791,66) (68921,66) (132651,36) (226981,43) (357911,46)};
\end{axis}
\end{tikzpicture}
\end{subfigure}
\hfill  \vspace*{.2cm} %
}\\%
\scalebox{1}{%
\hspace*{-.5cm}
\begin{subfigure}[t]{136.9pt}
\centering
\begin{tikzpicture}
\begin{axis}[
    xlabel={DOFs},
    ymode=log,
    width=\linewidth,
    height=4cm,
    grid=both,
    tick label style={font=\small},
    label style={font=\small},
    title style={font=\small},
    cycle list name=mycyclelist,
    every x tick scale label/.style={xshift = 2.9cm,yshift=-1.2cm}, %
]
\addplot coordinates {(1331,0.02) (9261,1.2) (29791,15) (68921,120)};
\addplot coordinates {(1331,0.033) (9261,2) (29791,26) (68921,85)};
\addplot coordinates {(1331,0.098) (9261,6.3) (29791,11) (68921,12) (132651,39) (226981,160) (357911,310)};
\addplot coordinates {(1331,0.082) (9261,1.5) (29791,5) (68921,9) (132651,10) (226981,40) (357911,99)};
\addplot coordinates {(1331,0.12) (9261,1) (29791,2.9) (68921,6.5) (132651,15) (226981,30) (357911,54)};
\addplot coordinates {(1331,0.18) (9261,1.6) (29791,4.4) (68921,11) (132651,25) (226981,50) (357911,89)};
\addplot coordinates {(1331,0.18) (9261,1.6) (29791,4.4) (68921,11) (132651,25) (226981,50) (357911,89)};
\end{axis}
\end{tikzpicture}
\end{subfigure}\hspace{-.6cm}
\begin{subfigure}[t]{136.9pt}
\centering
\begin{tikzpicture}
\begin{axis}[
    xlabel={DOFs},
    width=\linewidth,
    height=4cm,
    grid=both,
    tick label style={font=\small},
    label style={font=\small},
    title style={font=\small},
    cycle list name=mycyclelist,
    legend style={font=\tiny, at={(1,1.2)}, anchor=north east},
    legend columns=1,
    every x tick scale label/.style={xshift = 2.9cm,yshift=-.4cm}, %
]
\addplot coordinates {(1331,13) (9261,12) (29791,10) (68921,11)};
\addplot coordinates {(1331,13) (9261,12) (29791,10) (68921,10)};
\addplot coordinates {(1331,18) (9261,20) (29791,16) (68921,15) (132651,23) (226981,27) (357911,21)};
\addplot coordinates {(1331,21) (9261,17) (29791,13) (68921,16) (132651,14) (226981,16) (357911,19)};
\addplot coordinates {(1331,19) (9261,15) (29791,11) (68921,11) (132651,11) (226981,12) (357911,13)};
\addplot coordinates {(1331,19) (9261,15) (29791,11) (68921,11) (132651,11) (226981,12) (357911,13)};
\addplot coordinates {(1331,19) (9261,15) (29791,11) (68921,11) (132651,11) (226981,12) (357911,13)};
\end{axis}
\end{tikzpicture}
\end{subfigure}\hspace{-.9cm}
\begin{subfigure}[t]{136.9pt}
\centering
\begin{tikzpicture}
\begin{axis}[
    xlabel={DOFs},
    ymode=log,
    width=\linewidth,
    height=4cm,
    grid=both,
    tick label style={font=\small},
    label style={font=\small},
    title style={font=\small},
    cycle list name=mycyclelist,
    legend style={font=\tiny},
    legend pos=north west,
    every x tick scale label/.style={xshift = 2.9cm,yshift=1.65cm}, %
]
\addplot coordinates {};
\addplot coordinates {};
\addplot coordinates {(1331,1190) (9261,1537) (29791,1711) (68921,1097) (132651,2075) (226981,4283) (357911,4903)};
\addplot coordinates {(1331,636) (9261,587) (29791,1202) (68921,1121) (132651,709) (226981,1552) (357911,2221)};
\addplot coordinates {(1331,111) (9261,100) (29791,74) (68921,71) (132651,71) (226981,73) (357911,82)};
\addplot coordinates {(1331,112) (9261,101) (29791,74) (68921,71) (132651,73) (226981,74) (357911,82)};
\addplot coordinates {(1331,114) (9261,101) (29791,74) (68921,71) (132651,73) (226981,74) (357911,82)};
\end{axis}
\end{tikzpicture}
\end{subfigure}
\hfill  %
}%
\caption{Advection-Diffusion-Reaction equation: Solution of condensed system with outer tolerance $\mathrm{cgtol}=10^{-4}$ and $\lambda \in\{10^{-1},10^{-2},10^{-3}\}$ (top to bottom).}
\label{fig:condensed_ad}
\end{figure}

The numerical results for varying regularization parameters $\lambda$ and varying mesh sizes are shown in Figure~\ref{fig:condensed_ad}. Therein, a smaller $\lambda$ implies a larger condition number of the operator \eqref{eq:symmetricmatrix}, hence leading to a higher number of outer CG iterations.

\begin{figure}[htb]
    \centering

\begin{subfigure}[t]{136.9pt}
\centering
\begin{tikzpicture}
\begin{axis}[
    ylabel={Time [s]},
    xlabel={DOFs},
    width=\linewidth,
    height=4cm,
    grid=both,
    tick label style={font=\small},
    label style={font=\small},
    title style={font=\small},
    cycle list name=mycyclelist2,
    legend style={font=\small,at={(1.7,1.2)}, anchor=south},
    legend columns = 5,
    every x tick scale label/.style={xshift = 2.9cm,yshift=-.5cm}, %
]
\addplot coordinates {(1331,0.052) (9261,1.9) (29791,3.4) (68921,8.6) (132651,16) (226981,24) (357911,42)};
\addlegendentry{GMRES}
\addplot coordinates {(1331,0.045) (9261,0.86) (29791,2) (68921,3) (132651,6.2) (226981,13) (357911,21)};
\addlegendentry{GMRES(IC)}
\addplot coordinates {(1331,0.066) (9261,0.59) (29791,1.6) (68921,2.7) (132651,8.3) (226981,11) (357911,19)};
\addlegendentry{GMRES(MG)}
\addplot coordinates {(1331,0.1) (9261,0.84) (29791,2.6) (68921,4.1) (132651,13) (226981,17) (357911,30)};
\addlegendentry{Rapoport}
\addplot coordinates {(1331,0.095) (9261,0.84) (29791,2.6) (68921,4.1) (132651,13) (226981,17) (357911,30)};
\addlegendentry{Widlund}
\end{axis}
\end{tikzpicture}
\end{subfigure}\hspace{-.9cm}
\begin{subfigure}[t]{136.9pt}
\centering
\begin{tikzpicture}
\begin{axis}[
    xlabel={DOFs},
    width=\linewidth,
    height=4cm,
    grid=both,
    tick label style={font=\small},
    label style={font=\small},
    title style={font=\small},
    cycle list name=mycyclelist2,
    every x tick scale label/.style={xshift = 2.9cm,yshift=-.5cm}, %
]
\addplot coordinates {(1331,0.061) (9261,2.5) (29791,4.4) (68921,10) (132651,18) (226981,35) (357911,77)};
\addplot coordinates {(1331,0.054) (9261,0.98) (29791,2) (68921,5) (132651,6.5) (226981,22) (357911,71)};
\addplot coordinates {(1331,0.081) (9261,0.66) (29791,1.7) (68921,4.5) (132651,9.4) (226981,24) (357911,35)};
\addplot coordinates {(1331,0.12) (9261,0.98) (29791,2.8) (68921,7) (132651,15) (226981,38) (357911,52)};
\addplot coordinates {(1331,0.12) (9261,1) (29791,2.8) (68921,7) (132651,15) (226981,38) (357911,52)};
\end{axis}
\end{tikzpicture}
\end{subfigure}\hspace{-.9cm}
\begin{subfigure}[t]{136.9pt}
\centering
\begin{tikzpicture}
\begin{axis}[
    xlabel={DOFs},
    width=\linewidth,
    height=4cm,
    grid=both,
    tick label style={font=\small},
    label style={font=\small},
    title style={font=\small},
    cycle list name=mycyclelist2,
    every x tick scale label/.style={xshift = 2.9cm,yshift=-.5cm}, %
]
\addplot coordinates {(1331,0.098) (9261,6.3) (29791,11) (68921,12) (132651,39) (226981,160) (357911,310)};
\addplot coordinates {(1331,0.082) (9261,1.5) (29791,5) (68921,9) (132651,10) (226981,40) (357911,99)};
\addplot coordinates {(1331,0.12) (9261,1) (29791,2.9) (68921,6.5) (132651,15) (226981,30) (357911,54)};
\addplot coordinates {(1331,0.18) (9261,1.6) (29791,4.4) (68921,11) (132651,25) (226981,50) (357911,89)};
\addplot coordinates {(1331,0.18) (9261,1.6) (29791,4.4) (68921,11) (132651,25) (226981,50) (357911,89)};
\end{axis}
\end{tikzpicture}
\end{subfigure}\hspace{-.6cm}

     \caption{Time of the total outer iteration for indirect inner solvers. Outer tolerance $\mathrm{cgtol}=10^{-4}$ and $\lambda \in\{10^{-1},10^{-2},10^{-3}\}$ (left to right)}
    \label{fig:condensed_ad_time}
\end{figure}
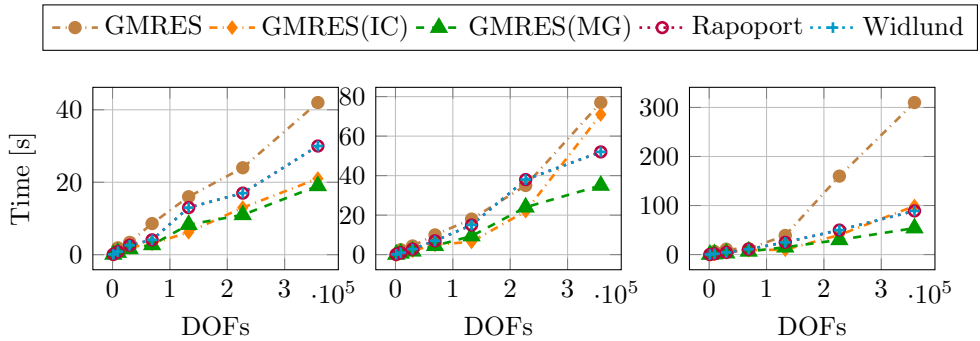

\noindent In the left column of Figure~\ref{fig:condensed_ad}, we depict the total time required for the full outer iteration. This includes the assembly of preconditioners and factorizations. We observe that the direct methods (LU and ILU factorizations) fail after already a few refinements of the grid as to be expected for the three-dimensional setting. The iterative inner solvers all perform similarly in this log scale, such that we provide a comparison with a linear scale in Figure~\ref{fig:condensed_ad_time}. Therein, observe that the multigrid-preconditioned GMRES method performs best (being an optimal method), while GMRES without preconditioner performs worst. Widlund's and Rapoport's method lead to very similar iteration numbers which is to be expected due to their strong similarity, see Section~\ref{sec:solvers}.

In the middle column of Figure~\ref{fig:condensed_ad}, we show the total number of outer CG iterations. We observe that GMRES and the incomplete-Cholesky-preconditioned GMRES lead to a higher number of outer iterations. We expect this to be due to an error estimator in the inner solves leading to termination before meeting the specified tolerance of $\mathrm{cgtol}/10$.

The right column of Figure~\ref{fig:condensed_ad} provides the total number of inner iterations for state and adjoint solves. Unpreconditioned GMRES has the highest iteration numbers (as to be expected), followed by the incomplete-Cholesky-preconditioned GMRES. The reason for this is that the incomplete Cholesky factorization with fixed fill-in does not represent a mesh-independent preconditioner (in contrast to the algebraic multigrid method), hence performing significantly worse for a larger number of spatial grid points.

Thus, we may conclude that GMRES, Rapoport's and Widlund's method preconditioned with an algebraic multigrid for the symmetric part yield the best inner solvers, both in terms of required iterations and total time. GMRES with incomplete-Cholesky preconditioner performs well in view of total required time, however requires significantly more inner iterations.

\subsubsection{Stokes equation}
As a second example, we consider the two-dimensional pressure-stabilized Stokes equation of Subsection~\ref{subsec:stokes}. We chose again full observation and a source term $f\equiv \begin{bmatrix}1 & 1\end{bmatrix}^\top$. We use the same parameters as in Table~\ref{tab:methods_params} with the difference that we now use four multigrid V-cycles and an incomplete Cholesky with drop tolerance $10^{-2}$.

\begin{figure}[htb]
\centering
\begin{subfigure}[t]{136.9pt}
\centering
\begin{tikzpicture}
\begin{axis}[    
    title={Time [s]},
    xlabel={DOFs},
    ymode=log,
    width=\linewidth,
    height=4cm,
    grid=both,
    tick label style={font=\small},
    label style={font=\small},
    title style={font=\small},
    cycle list name=mycyclelist,
    legend style={font=\small,at={(1.7,1.35)}, anchor=south},
    legend columns = 4,
    every x tick scale label/.style={xshift = 3.4cm,yshift=-1.13cm}, %
]
\addplot coordinates {(1202,0.035) (4182,0.62) (8962,4.5) (15542,24) (23922,83)};
\addlegendentry{Direct}
\addplot coordinates {(1202,0.061) (4182,1.2) (8962,7) (15542,24) (23922,67)};
\addlegendentry{ILU}
\addplot coordinates {(1202,0.12) (4182,7.1) (8962,25) (15542,66) (23922,45) (34102,80) (46082,120) (59862,150) (92822,280)};
\addlegendentry{GMRES}
\addplot coordinates {(1202,0.04) (4182,0.41) (8962,1.1) (15542,2.3) (23922,2.2) (34102,3.7) (46082,4.4) (59862,6.7) (92822,11)};
\addlegendentry{GMRES(IC)}
\addplot coordinates {(1202,0.1) (4182,0.36) (8962,0.65) (15542,1.1) (23922,1.5) (34102,2.1) (46082,2.4) (59862,2.8) (92822,4.2)};
\addlegendentry{GMRES(MG)}
\addplot coordinates {(1202,0.18) (4182,0.37) (8962,1.3) (15542,1) (23922,3.7) (34102,5.2) (46082,6.6) (59862,8.3)};
\addlegendentry{Rapoport}
\addplot coordinates {(1202,0.18) (4182,0.38) (8962,1.3) (15542,1) (23922,3.6) (34102,5.2) (46082,6.6) (59862,8.3)};
\addlegendentry{Widlund}
\end{axis}
\end{tikzpicture}
\end{subfigure}\hspace{-.9cm}
\begin{subfigure}[t]{136.9pt}
\centering
\begin{tikzpicture}
\begin{axis}[
    title={Outer CG iterations},
    xlabel={DOFs},
    width=\linewidth,
    height=4cm,
    grid=both,
    tick label style={font=\small},
    label style={font=\small},
    title style={font=\small},
    cycle list name=mycyclelist,
    every x tick scale label/.style={xshift = 3.4cm,yshift=-.37cm}, %
]
\addplot coordinates {(1202,9) (4182,8) (8962,7) (15542,7) (23922,6)};
\addplot coordinates {(1202,10) (4182,9) (8962,9) (15542,8) (23922,7)};
\addplot coordinates {(1202,10) (4182,9) (8962,8) (15542,8) (23922,6) (34102,6) (46082,6) (59862,5) (92822,5)};
\addplot coordinates {(1202,10) (4182,9) (8962,8) (15542,7) (23922,7) (34102,7) (46082,6) (59862,6) (92822,6)};
\addplot coordinates {(1202,10) (4182,9) (8962,8) (15542,8) (23922,7) (34102,7) (46082,6) (59862,6) (92822,5)};
\addplot coordinates {(1202,11) (4182,8) (8962,7) (15542,7) (23922,11) (34102,11) (46082,9) (59862,9)};
\addplot coordinates {(1202,11) (4182,8) (8962,7) (15542,7) (23922,11) (34102,11) (46082,9) (59862,9)};
\end{axis}
\end{tikzpicture}
\end{subfigure}\hspace{-.9cm}
\begin{subfigure}[t]{136.9pt}
\centering
\begin{tikzpicture}
\begin{axis}[
    title={Inner iterations},
    xlabel={DOFs},
    ymode=log,
    width=\linewidth,
    height=4cm,
    grid=both,
    tick label style={font=\small},
    label style={font=\small},
    title style={font=\small},
    cycle list name=mycyclelist,
    every x tick scale label/.style={xshift = 3.4cm,yshift=1.33cm}, %
]
\addplot coordinates {};
\addplot coordinates {};
\addplot coordinates {(1202,1148) (4182,1916) (8962,2424) (15542,3122) (23922,3022) (34102,3612) (46082,4161) (59862,4070) (92822,5019)};
\addplot coordinates {(1202,218) (4182,322) (8962,412) (15542,472) (23922,542) (34102,613) (46082,598) (59862,681) (92822,774)};
\addplot coordinates {(1202,104) (4182,94) (8962,84) (15542,84) (23922,74) (34102,75) (46082,66) (59862,56) (92822,56)};
\addplot coordinates {(1202,112) (4182,61) (8962,115) (15542,49) (23922,117) (34102,117) (46082,110) (59862,106)};
\addplot coordinates {(1202,112) (4182,61) (8962,116) (15542,49) (23922,117) (34102,117) (46082,110) (59862,106)};
\end{axis}
\end{tikzpicture}
\end{subfigure}

\vspace*{.2cm}
\begin{subfigure}[t]{136.9pt}
\centering
\begin{tikzpicture}
\begin{axis}[
    xlabel={$n$},
    ymode=log,
    width=\linewidth,
    height=4cm,
    grid=both,
    tick label style={font=\small},
    label style={font=\small},
    title style={font=\small},
    cycle list name=mycyclelist,
    every x tick scale label/.style={xshift = 3.4cm,yshift=-1.13cm}
]
\addplot coordinates {(1202,0.035) (4182,0.68) (8962,4.5) (15542,24) (23922,83) (34102,230)};
\addplot coordinates {(1202,0.077) (4182,0.85)};
\addplot coordinates {(1202,0.14) (4182,6.7) (8962,30) (15542,82) (23922,69) (34102,100) (46082,150) (59862,220) (92822,320)};
\addplot coordinates {(1202,0.04) (4182,0.4) (8962,1.1) (15542,3.2) (23922,3.3) (34102,4.2) (46082,6.8) (59862,9.9) (92822,13)};
\addplot coordinates {(1202,0.11) (4182,0.37) (8962,0.69) (15542,1.2) (23922,1.5) (34102,2.3) (46082,2.7) (59862,2.9) (92822,3.9)};
\addplot coordinates {(1202,0.19) (4182, 0.68) (15542, 0.96)};
\addplot coordinates {(1202,0.18) (4182, 0.67) (15542, 0.96)};
\end{axis}
\end{tikzpicture}
\end{subfigure}\hspace{-.9cm}
\begin{subfigure}[t]{136.9pt}
\centering
\begin{tikzpicture}
\begin{axis}[
    xlabel={DOFs},
    width=\linewidth,
    height=4cm,
    grid=both,
    tick label style={font=\small},
    label style={font=\small},
    title style={font=\small},
    cycle list name=mycyclelist,
    legend style={font=\tiny, at={(1,1.1)}, anchor=north east},
    legend columns=1,
    every x tick scale label/.style={xshift = 3.4cm,yshift=-.37cm}
]
\addplot coordinates {(1202,9) (4182,8) (8962,7) (15542,7) (23922,6) (34102,6)};
\addplot coordinates {(1202,10) (4182,9)};
\addplot coordinates {(1202,10) (4182,8) (8962,7) (15542,7) (23922,6) (34102,6) (46082,5) (59862,5) (92822,4)};
\addplot coordinates {(1202,10) (4182,9) (8962,8) (15542,7) (23922,7) (34102,6) (46082,6) (59862,6) (92822,4)};
\addplot coordinates {(1202,10) (4182,9) (8962,8) (15542,7) (23922,6) (34102,6) (46082,6) (59862,5) (92822,4)};
\addplot coordinates {(1202,11) (4182,7) (15542,6)};
\addplot coordinates {(1202,11) (4182,7) (15542,6)};
\end{axis}
\end{tikzpicture}
\end{subfigure}\hspace{-.9cm}
\begin{subfigure}[t]{136.9pt}
\centering
\begin{tikzpicture}
\begin{axis}[
    xlabel={DOFs},
    ymode=log,
    width=\linewidth,
    height=4cm,
    grid=both,
    tick label style={font=\small},
    label style={font=\small},
    title style={font=\small},
    cycle list name=mycyclelist,
    every x tick scale label/.style={xshift = 3.4cm,yshift=1.33cm},
]
\addplot coordinates {};
\addplot coordinates {};
\addplot coordinates {(1202,1300) (4182,1958) (8962,2494) (15542,3261) (23922,3537) (34102,4211) (46082,4132) (59862,4698) (92822,4752)};
\addplot coordinates {(1202,245) (4182,375) (8962,472) (15542,538) (23922,658) (34102,679) (46082,776) (59862,870) (92822,730)};
\addplot coordinates {(1202,115) (4182,106) (8962,95) (15542,87) (23922,76) (34102,85) (46082,74) (59862,65) (92822,54)};
\addplot coordinates {(1202,126) (4182,137) (15542,48)};
\addplot coordinates {(1202,126) (4182,137) (15542,48)};
\end{axis}
\end{tikzpicture}
\end{subfigure} \caption{Stabilized Stokes equation: Solution of condensed system with outer tolerance $10^{-4}$ and $\lambda =10^{-1}$. Above: $s_1=s_2=1$. Below: $s_1=0$, $s_2=1$. }
\label{fig:condensed_stokes}
\end{figure}
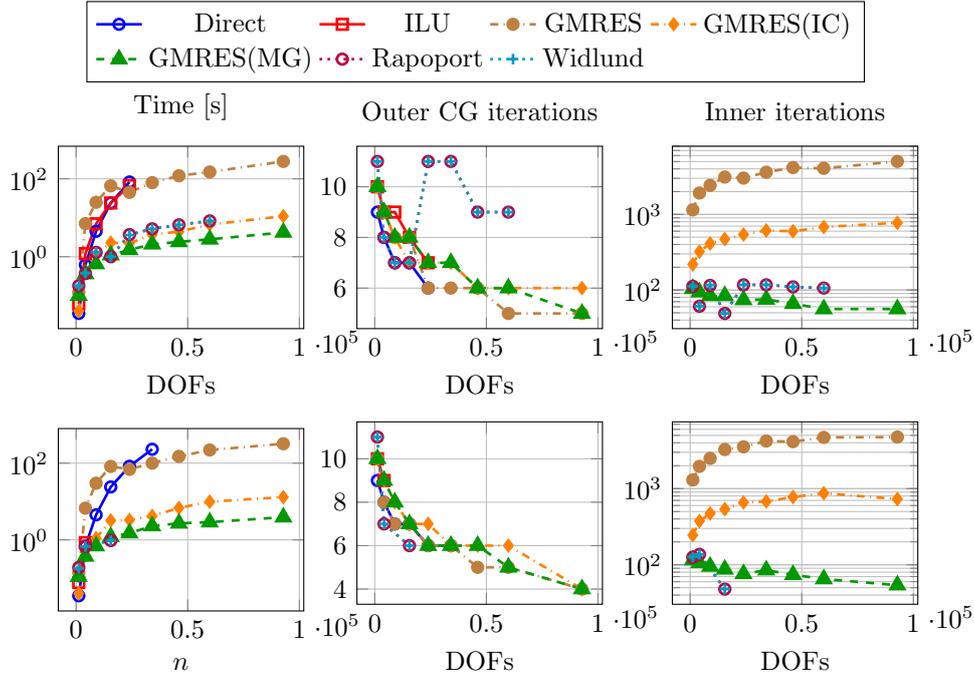

The results are depicted in Figure~\ref{fig:condensed_stokes} for the stabilized case. The upper row depicts the case of pressure regularization in the full $H^1$ norm, while the lower row only includes the gradient term. Again, we observe that the direct methods fail after a few refinements. Unpreconditioned GMRES and GMRES preconditioned by an incomplete Cholesky factorization performs well in terms of outer CG iterations, however requires a high amount of computation time due to the high number of inner iterations. Again, GMRES endowed with an AMG preconditioner for the symmetric part performs best. Rapoport's and Widlund's method fail due to lack of symmetry after a few refinement steps due to the inexact evaluation of the preconditioner by the multigrid method, motivating future research using flexible methods as suggested in~\cite{diab2023flexible}.

\subsection{Projected conjugate gradients}\label{sec:ppcg}
As another method, we suggest a conjugate gradient variant applied to the full optimality system~\eqref{eq:optsys} without eliminating the state variable via the state-control map. 

As the governing operator is not positive definite (being a saddle point matrix and due to the zero in the bottom right block), we use the constraint preconditioner
\begin{equation*}
        \mathcal{P} := \begin{bmatrix}
        0 & 0 & \calA^*\\
        0 & \lambda I & -\calB^*\\
        \calA& -\calB &0
        \end{bmatrix}
\end{equation*}
which ensures that the iterates $\begin{bmatrix}x_k &u_k\end{bmatrix}^\top$ evolve in the kernel of the equality constraint.  Thus, the proposed method falls under the class of projected preconditioned conjugate gradient (PPCG) methods. Such a preconditioner was used e.g.\ in \cite{lubkoll2017affine} in the context of affine covariant composite step methods. This preconditioner shares two important features: \\
1.    If the conjugate gradient iteration is started with $\begin{bmatrix}x^0&u^0&p^0\end{bmatrix}^\top$ such that the primal variables are admissible, i.e., that $\calA x^0 + \calB u^0 = 0$, all residuals vanish in the last block component such that the preconditioned system gives rise to search directions $\begin{bmatrix}d_x&d_u&d_p\end{bmatrix}^\top$ with $\calA d_x+\calB d_u = 0$. Hence, the iteration yields a gradient method fully evolving in the kernel of the equality constraint, hence removing indefiniteness.
    \\
    2. If a solver for $\mathcal{A}$ and $\mathcal{A^*}$ is available, then this preconditioner can be efficiently evaluated in a block-row-wise fashion starting with the top left block, i.e., by the preconditioner
\[
\calP^{-1}
=
\begin{bmatrix}
\frac{1}{\lambda} \calA^{-1} \calB \calB^* \calA^{-*} & \frac{1}{\lambda} \calA^{-1} \calB & \calA^{-1} \\
\frac{1}{\lambda} \calB^*\calA^{-*} & \frac{1}{\lambda} I & 0 \\
\calA^{-*} & 0 & 0
\end{bmatrix}.
\]
In this way  the preconditioned version of the optimality system \eqref{eq:optsys} is given by
\begin{equation*}
    \calP^{-1} \begin{bmatrix}
    \calC^*\calC & 0 & \calA^*\\
    0 & \lambda I & -\calB^*\\
    \calA& -\calB &0
\end{bmatrix} = 
\begin{bmatrix}
    \frac{1}{\lambda} \calA^{-1} \calB \calB^* \calA^{-*}\calC^* \calC + I & 0 & 0\\
    \frac{1}{\lambda} \calB^*\calA^{-1} \calC^*\calC & I & 0\\
    \calA^{-1}\calC^*\calC & 0 & I
\end{bmatrix}
\end{equation*}
is a bounded block operator as it does not involve $\calA$ but only its inverses, which in case of differential equations have a smoothing effect. Further, the efficient numerical evaluation of this preconditioner requires only one solve with $\calA$ and  $\calA^{*}$.
  
In view of the point above, we will again use the solvers tailored to the structure of $\mathcal{A} = \calH + \calS$ such as (preconditioned) GMRES, Widlund's or Rapoport's method. As these approaches constitute iterative methods, inexact evaluations of the preconditioner are inevitable. In this work, to focus on the main contribution, we will set a very small solver tolerance, such that the inner solves are performed with high accuracy, see the setting described in Table~\ref{tab:methods_ppcg}. However,  various works have considered inexact CG methods. A general framework for inexact evaluations of the preconditioner building upon a bound on the corresponding residual was developed in \cite{GoluYe99}. %
Moreover, inexact PPCG methods in the context of inexact SQP methods were discussed in \cite{HeinRidz14}, and in \cite{Schiela2021} a primal-dual projection method to alleviate inexact solves is presented. 

\begin{table}[htb]
\centering
\scalebox{.75}{
\begin{tabular}{c | c c c c c c c}
\hline
& direct & ILU & GMRES & GMRES(IC)  & GMRES(MG)  & Rapoport & Widlund  \\
\hline
assembly & lu(A) & ilu(A) & - &  ichol(H,$10^{-1}$)& multigrid (2c) & multigrid & multigrid  \\
tolerance & - & cgtol/10 (drop) & 1e-6 & 1e-6 & 1e-6  & 1e-6 &1e-6   \\
\hline
\end{tabular}
}
\caption{Used inner solver for evaluation of the constraint preconditioner for the outer CG iteration with tolerance $\mathrm{cgtol}$.}
\label{tab:methods_ppcg}
\end{table}

\begin{figure}[htbp]
\centering
\begin{subfigure}[t]{136.9pt}
\centering
\begin{tikzpicture}
\begin{axis}[
    xlabel={DOFs},
    title={Time [s]},
    ymode=log,
    width=\linewidth,
    height=4cm,
    grid=both,
    tick label style={font=\small},
    label style={font=\small},
    title style={font=\small},
    cycle list name=mycyclelist,
    legend style={font=\small,at={(1.7,1.35)}, anchor=south},
    legend columns = 4,
    every x tick scale label/.style={xshift = 3.1cm,yshift=-.8cm}, %
]
\addplot coordinates {(1331,0.15) (9261,2.6) (29791,20) (68921,130)};
\addlegendentry{Direct}
\addplot coordinates {(1331,0.16)};
\addlegendentry{ILU}
\addplot coordinates {(1331,0.21) (9261,9.4) (29791,24) (68921,53) (132651,140) (226981,210) (357911,540) (531441,1100)};
\addlegendentry{GMRES}
\addplot coordinates {(1331,0.18) (9261,3) (29791,11) (68921,30) (132651,61) (226981,140) (357911,240) (531441,590)};
\addlegendentry{GMRES(IC)}
\addplot coordinates {(1331,0.25) (9261,2.4) (29791,8.5) (68921,23) (132651,61) (226981,120) (357911,240) (531441,440)};
\addlegendentry{GMRES(MG)}
\addplot coordinates {(1331,0.32) (9261,2.9) (29791,11) (68921,28) (132651,73) (226981,150) (357911,270) (531441,520)};
\addlegendentry{Rapoport}
\addplot coordinates {(1331,0.32) (9261,2.9) (29791,11) (68921,27) (132651,71) (226981,150) (357911,270) (531441,510)};
\addlegendentry{Widlund}
\end{axis}
\end{tikzpicture}
\end{subfigure}\hspace{-.9cm}
\begin{subfigure}[t]{136.9pt}
\centering
\begin{tikzpicture}
\begin{axis}[
    xlabel={DOFs},
    title = {Outer CG iterations},
    width=\linewidth,
    height=4cm,
    grid=both,
    tick label style={font=\small},
    label style={font=\small},
    title style={font=\small},
    cycle list name=mycyclelist,
    every x tick scale label/.style={xshift = 3.1cm,yshift=-.38cm}, %
]
\addplot coordinates {(1331,9) (9261,9) (29791,7) (68921,9)};
\addplot coordinates {(1331,9)};
\addplot coordinates {(1331,11) (9261,12) (29791,11) (68921,9) (132651,11) (226981,9) (357911,11) (531441,11)};
\addplot coordinates {(1331,10) (9261,10) (29791,10) (68921,10) (132651,8) (226981,10) (357911,8) (531441,10)};
\addplot coordinates {(1331,10) (9261,9) (29791,8) (68921,8) (132651,8) (226981,8) (357911,8) (531441,8)};
\addplot coordinates {(1331,10) (9261,9) (29791,8) (68921,8) (132651,8) (226981,8) (357911,8) (531441,8)};
\addplot coordinates {(1331,10) (9261,9) (29791,8) (68921,8) (132651,8) (226981,8) (357911,8) (531441,8)};
\end{axis}
\end{tikzpicture}
\end{subfigure}\hspace{-.9cm}
\begin{subfigure}[t]{136.9pt}
\centering
\begin{tikzpicture}
\begin{axis}[
    xlabel={DOFs},
    ymode=log,
    title = {Inner iterations},
    width=\linewidth,
    height=4cm,
    grid=both,
    tick label style={font=\small},
    label style={font=\small},
    title style={font=\small},
    cycle list name=mycyclelist,
    legend style={font=\tiny},
    legend pos=north west,
    every x tick scale label/.style={xshift = 3.1cm,yshift=1.6cm}, %
]
\addplot coordinates {};
\addplot coordinates {};
\addplot coordinates {(1331,652) (9261,1497) (29791,2037) (68921,2106) (132651,3050) (226981,2908) (357911,4092) (531441,4503)};
\addplot coordinates {(1331,363) (9261,621) (29791,902) (68921,1152) (132651,1145) (226981,1592) (357911,1551) (531441,2199)};
\addplot coordinates {(1331,78) (9261,70) (29791,63) (68921,64) (132651,63) (226981,64) (357911,64) (531441,64)};
\addplot coordinates {(1331,78) (9261,70) (29791,64) (68921,63) (132651,63) (226981,63) (357911,64) (531441,64)};
\addplot coordinates {(1331,78) (9261,70) (29791,64) (68921,63) (132651,63) (226981,63) (357911,64) (531441,64)};
\end{axis}
\end{tikzpicture}
\end{subfigure}
\hspace*{-.08cm}
\begin{subfigure}[t]{136.9pt}
\centering
\begin{tikzpicture}
\begin{axis}[
    xlabel={DOFs},
    ymode=log,
    width=\linewidth,
    height=4cm,
    grid=both,
    tick label style={font=\small},
    label style={font=\small},
    title style={font=\small},
    cycle list name=mycyclelist,
    every x tick scale label/.style={xshift = 3.3cm,yshift=-1.02cm}, %
]
\addplot coordinates {(1202,0.046) (4182,0.66) (8962,6) (15542,27) (23922,86)};
\addplot coordinates {};
\addplot coordinates {(1202,0.16) (4182,13) (8962,59) (15542,130) (23922,140) (34102,230) (46082,370) (59862,540) (75442,600) (92822,950)};
\addplot coordinates {(1202,0.078) (4182,2.9) (8962,13) (15542,28) (23922,33) (34102,53) (46082,80) (59862,110) (75442,130) (92822,190)};
\addplot coordinates {(1202,0.11) (4182,0.43) (8962,2.3) (15542,4) (23922,6.2) (34102,9.1) (46082,12) (59862,16) (75442,21) (92822,26)};
\addplot coordinates {(1202,0.19) (4182,0.6) (8962,3.1) (15542,5.1) (23922,8.1) (34102,12) (46082,18) (59862,19) (75442,30) (92822,37)};
\addplot coordinates {(1202,0.18) (4182,0.59) (8962,3.1) (15542,5.1) (23922,8.1) (34102,12) (46082,18) (59862,19) (75442,30) (92822,37)};
\end{axis}
\end{tikzpicture}
\end{subfigure}\hspace{-.83cm}
\begin{subfigure}[t]{136.9pt}
\centering
\begin{tikzpicture}
\begin{axis}[
    xlabel={DOFs},
    width=\linewidth,
    height=4cm,
    grid=both,
    tick label style={font=\small},
    label style={font=\small},
    title style={font=\small},
    cycle list name=mycyclelist,
    legend style={font=\tiny, at={(1,1.2)}, anchor=north east},
    legend columns=1,
    every x tick scale label/.style={xshift = 3.3cm,yshift=-.39cm}, %
]
\addplot coordinates {(1202,9) (4182,9) (8962,9) (15542,9) (23922,9) (34102,9)};
\addplot coordinates {};
\addplot coordinates {(1202,9) (4182,9) (8962,9) (15542,9) (23922,9) (34102,9) (46082,9) (59862,9) (75442,9) (92822,9)};
\addplot coordinates {(1202,10) (4182,10) (8962,10) (15542,10) (23922,10) (34102,10) (46082,10) (59862,10) (75442,10) (92822,10)};
\addplot coordinates {(1202,10) (4182,10) (8962,9) (15542,9) (23922,9) (34102,9) (46082,9) (59862,9) (75442,9) (92822,9)};
\addplot coordinates {(1202,11) (4182,10) (8962,11) (15542,10) (23922,10) (34102,10) (46082,11) (59862,9) (75442,11) (92822,11)};
\addplot coordinates {(1202,11) (4182,10) (8962,11) (15542,10) (23922,10) (34102,10) (46082,11) (59862,9) (75442,11) (92822,11.0)};
\end{axis}
\end{tikzpicture}
\end{subfigure}\hspace{-.9cm}
\begin{subfigure}[t]{136.9pt}
\centering
\begin{tikzpicture}
\begin{axis}[
    xlabel={DOFs},
    ymode=log,
    width=\linewidth,
    height=4cm,
    grid=both,
    tick label style={font=\small},
    label style={font=\small},
    title style={font=\small},
    cycle list name=mycyclelist,
    legend style={font=\tiny},
    legend pos=north west,
    every x tick scale label/.style={xshift = 3.3cm,yshift=1.585cm}, %
]
\addplot coordinates {};
\addplot coordinates {};
\addplot coordinates {(1202,1326) (4182,2425) (8962,3528) (15542,4685) (23922,5795) (34102,6917) (46082,7996) (59862,9040) (75442,10059) (92822,11092)};
\addplot coordinates {(1202,716) (4182,1282) (8962,1856) (15542,2408) (23922,2973) (34102,3509) (46082,4031) (59862,4571) (75442,5115) (92822,5619)};
\addplot coordinates {(1202,115) (4182,114) (8962,102) (15542,103) (23922,103) (34102,102) (46082,101) (59862,101) (75442,102) (92822,101)};
\addplot coordinates {(1202,126) (4182,116) (8962,126) (15542,114) (23922,114) (34102,116) (46082,127) (59862,106) (75442,127) (92822,126)};
\addplot coordinates {(1202,125) (4182,116) (8962,126) (15542,114) (23922,114) (34102,116) (46082,127) (59862,106) (75442,127) (92822,126)};
\end{axis}
\end{tikzpicture}
\end{subfigure} \caption{PPCG: Solution of optimality system with PPCG for advection-diffusion-reaction equation with outer tolerance $10^{-4}$ and $\lambda =10^{-4}$ (above) and Stokes equation with outer tolerance $10^{-3}$ and $\lambda = 10^{-5}$ (below).}
\label{fig:ppcg}
\end{figure}
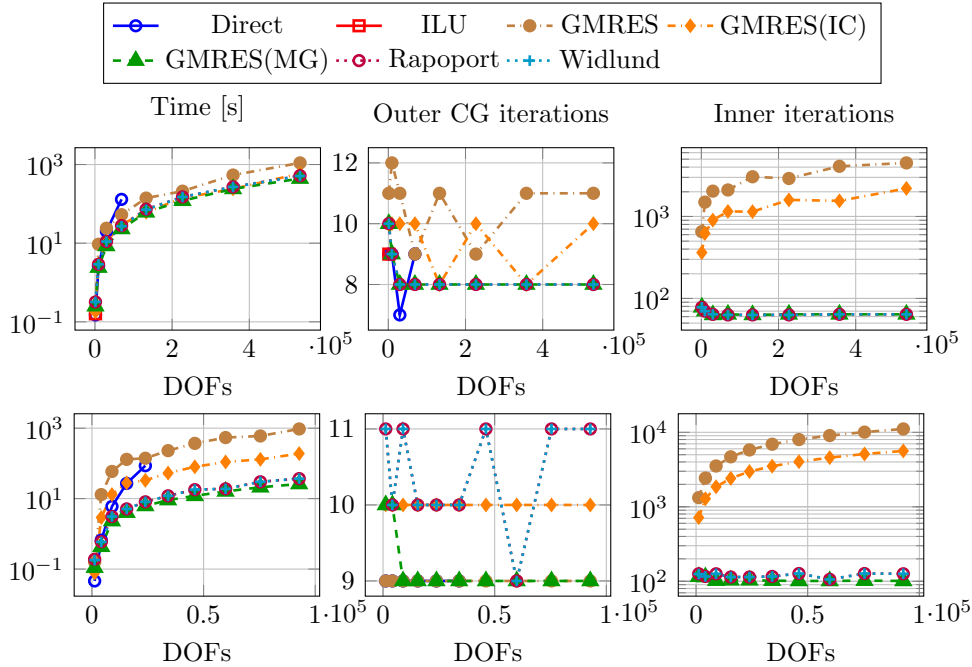

\noindent In Figure~\ref{fig:ppcg}, we depict the resulting runtime and iteration numbers for the outer CG loop applied to the optimality system \eqref{eq:optsys} endowed with the proposed constraint preconditioner. 
\noindent In the top row of Figure~\ref{fig:ppcg}, we show the results for the advection-diffusion-reaction equation. We see that the total number of outer iterations is robust w.r.t.\ refinements. The reason for this is that the constraint preconditioned operator is bounded and hence, the discretization has a uniformly bounded condition number. However, as in the approach with state elimination suggested in Subsection~\ref{sec:reduced}, the time and inner iterations required by GMRES and the incomplete-Cholesky-preconditioned variant is significantly higher as for the other methods. Thus, we conclude that also in this approach using a constraint preconditioner, again the tailored methods being GMRES, Rapoport's methods or Widlund's method endowed with an algebraic multigrid preconditioner for the symmetric part perform best.

\section{Conclusion}
We have analyzed and evaluated tailored methods and preconditioning for accretive systems occurring in partial differential equations such as dissipative and port-Hamiltonian systems. Therein, the preconditioner is obtained by the symmetric part of the underlying operator, allowing for short-recurrence Krylov subspace methods using suitable inner products such as Widlund's and Rapoport's method. In the first part of this work, we have analyzed this approach in view of operator preconditioning for a wide range of differential operators occurring in advection-diffusion-reaction equations, fluid mechanics, wave propagation or elasticity. In the second part we have suggested two approaches to include these methods in a large scale optimal control solver, i.e., a reduced formulation and projected preconditioned cg method for the full optimality system. We have illustrated that approximating the symmetric preconditioner via algebraic multigrid methods or incomplete Cholesky factorization leads to highly efficient and robust solvers. As future work we plan to analyze this approach for nonlinear partial differential operators and abstract differential-algebraic operators.

\bibliographystyle{abbrv}
\bibliography{refs.bib}

\end{document}